\documentclass[12pt]{article}
\usepackage{amsmath}
\usepackage{amssymb}
\usepackage{vatola}

\def\q{\quad}
\def\qq{\qquad}
\def\qtq#1{\q\t{#1}\q}
\def\mod#1{\ (\text{\rm mod}\ #1)}
\def\t{\text}
\def\f{\frac}
\def\e{\equiv}
\def\b{\binom}

\def\sls#1#2{(\f{#1}{#2})}
 \def\ls#1#2{\big(\f{#1}{#2}\big)}
\def\Ls#1#2{\Big(\f{#1}{#2}\Big)}

\let \pro=\proclaim
\let \endpro=\endproclaim

\begin{document}
 \centerline {\bf
Conjectures on congruences involving binomial coefficients and
Ap\'ery-like numbers}
\par\q\newline
\centerline{Zhi-Hong Sun}\newline \centerline{School of Mathematics
and Statistics}\centerline{Huaiyin Normal University}
\centerline{Huaian, Jiangsu 223300, P.R. China} \centerline{Email:
zhsun@hytc.edu.cn} \centerline{Homepage:
http://maths.hytc.edu.cn/szh1.htm}
 \abstract{In this paper, we pose many challenging conjectures on
 congruences involving binomial coefficients and Ap\'ery-like numbers.
 \par\q
 }
 \endabstract
\section*{1. Introduction}
\par\q
\par  We first illustrate the notation in the paper. Let $\Bbb Z$ be the set of integers, and for a
prime $p$ let $\Bbb Z_p$ be the set of rational numbers whose
denominator is not divisible by $p$. Let $[x]$ be
 the greatest integer not exceeding $x$. For $a\in\Bbb Z$ and given
odd prime $p$ let $\sls ap$ denote the Legendre symbol.  For
positive integers $a,b$ and $n$, if $n=ax^2+by^2$ for some integers
$x$ and $y$, we briefly write that $n=ax^2+by^2$.
\par Let $p$ be an odd prime. It is clear that
$$\b{2k}k=\f{2k(2k-1)\cdots(k+1)}{k!}\e 0\mod p\qtq{for}k=\f{p+1}2,\ldots,p-1$$
and so $\f 1{2k-1}\b{2k}k\in\Bbb Z_p$ for $k=0,1,\ldots,p-1$.
Actually, for $k\ge 1$,
$$\f 1{2k-1}\b{2k}k=2\Big(\b{2k-2}{k-1}-\b{2k-2}k\Big)
=\f 2k\b{2k-2}{k-1}=2C_{k-1}\in\Bbb Z,$$
where $C_k=\f 1{k+1}\b{2k}k$ is the $k$-th Catalan number.

\par
Let $p>3$ be a prime. In 1987, Beukers[B] conjectured a congruence
equivalent to
$$\sum_{k=0}^{p-1}\f{\b{2k}k^3}{64^k}\e
\cases 0\mod{p^2}&\t{if $p\e 3\mod 4$,}\\4x^2-2p\mod{p^2}&\t{if
$p=x^2+4y^2\e 1\mod 4$.}\endcases$$ This congruence was proved by
several authors including Ishikawa[I]($p\e 1\mod 4$), van
Hamme[vH]($p\e 3\mod 4$) and Ahlgren[A]. In 2003,
Rodriguez-Villegas[RV] posed 22 conjectures on supercongruences
modulo $p^2$. In particular, the following congruences are
equivalent to conjectures due to Rodriguez-Villegas:
$$\align &\sum_{k=0}^{p-1}\f{\b{2k}k^2\b{3k}k}{108^k}\e
\cases 4x^2-2p\mod{p^2}&\t{if $p=x^2+3y^2\e 1\mod 3$,}
\\0\mod{p^2}&\t{if $p\e 2\mod 3$,}
\endcases
\\&\sum_{k=0}^{p-1}\f{\b{2k}k^2\b{4k}{2k}}{256^k}\e
\cases 4x^2-2p\mod{p^2}&\t{if $p=x^2+2y^2\e 1,3\mod 8$,}
\\0\mod{p^2}&\t{if $p\e 5,7\mod 8$,}
\endcases
\\&\Ls p3\sum_{k=0}^{p-1}\f{\b{2k}k\b{3k}k\b{6k}{3k}}{12^{3k}}\e
\cases 4x^2-2p\mod{p^2}&\t{if $p=x^2+4y^2\e 1\mod 4$,}
\\0\mod{p^2}&\t{if $p\e 3\mod 4$.}
\endcases\endalign$$
  These conjectures have been solved  by
Mortenson[M] and Zhi-Wei Sun[Su2]. Since
$$C_k=\f{\b{2k}k}{k+1}=\b{2k}k-\b{2k}{k+1},$$
from [Su2] one may deduce that
$$\align
&\sum_{k=0}^{p-1}\f{\b{2k}k^3}{64^k(k+1)}\e \cases 4x^2-2p \mod
{p^2}\qq\t{if $p=x^2+4y^2\e 1\mod 4$,}
\\-(2p+2-2^{p-1})\b{(p-1)/2}{(p+1)/4}^2\mod {p^2}\q\t{if $p\e 3\mod
4$,}\endcases
\\&\sum_{k=0}^{p-1}\f{\b{2k}k^3}{(-8)^k(k+1)}\e 6x^2-4p \mod
{p^2}\qtq{for}p=x^2+4y^2\e 1\mod 4,
\\&\sum_{k=0}^{p-1}\f{\b{2k}k^2\b{3k}k}{108^k(k+1)}\e
 4x^2-2p\mod{p^2}\qtq{for}p=x^2+3y^2\e 1\mod 3,
\\&\sum_{k=0}^{p-1}\f{\b{2k}k^2\b{4k}{2k}}{256^k(k+1)}\e
 4x^2-2p\mod{p^2}\qtq{for}p=x^2+2y^2\e 1,3\mod 8,
\\&\Ls p3\sum_{k=0}^{p-1}\f{\b{2k}k\b{3k}k\b{6k}{3k}}{12^{3k}(k+1)}\e
 4x^2-2p\mod{p^2}\qtq{for}p=x^2+4y^2\e 1\mod 4.
\endalign$$

Let $p$ be an odd prime, $m\in\Bbb Z$ and $p\nmid m$. In [Su1,Su3],
  Z.W. Sun posed many conjectures concerning congruences modulo $p^2$
involving the sums
$$\sum_{k=0}^{p-1}\f{\b{2k}k^3}{m^k},\q
\sum_{k=0}^{p-1}\f{\b{2k}k^2\b{3k}k}{m^k},\q
\sum_{k=0}^{p-1}\f{\b{2k}k^2\b{4k}{2k}}{m^k}, \q
\sum_{k=0}^{p-1}\f{\b{2k}k\b{3k}k\b{6k}{3k}}{m^k}.$$ For 13 similar
conjectures see [S1].  Most of these congruences modulo $p$ were
proved by the author in [S2-S5]. In [S10], the author conjectured
many congruences modulo $p^3$ involving the above sums. For
instance, for any prime $p\not=2,3,7$,
$$\sum_{k=0}^{p-1}\b{2k}k^3\e
\cases  4x^2-2p-\f{p^2}{4x^2}\mod {p^3}\ \t{if $p\e 1,2,4\mod 7$ and
so $p=x^2+7y^2$,}
\\-\f {11}4p^2\b{3[p/7]}{[p/7]}^{-2} \e
-11p^2\b{[3p/7]}{[p/7]}^{-2}\mod {p^3}\qq\ \t{if $7\mid p-3$,}
\\-\f{99}{64}p^2\b{3[p/7]}{[p/7]}^{-2}\e
-11p^2\b{[6p/7]}{[2p/7]}^{-2}\mod {p^3}\qq\ \t{if $7\mid p-5$,}
\\-\f{25}{176}p^2\b{3[p/7]}{[p/7]}^{-2}\e
-11p^2\b{[3p/7]}{[p/7]+1}^{-2}\mod {p^3}\q\t{if $7\mid p-6$.}
\endcases$$
If we replace $\b{3[p/7]}{[p/7]}$ with $\b{[3p/7]}{[p/7]}$, for $p\e
3,5,6\mod 7$ we have
$$\sum_{k=0}^{p-1}\b{2k}k^3\e
\cases -\f {11p^2}{\b{[3p/7]}{[p/7]}^2}\mod {p^3} &\t{if $p\e 3\mod
7$,}
\\-\f {11p^2}{16\b{[3p/7]}{[p/7]}^2}\mod {p^3} &\t{if $p\e 5\mod
7$,}
\\-\f {11p^2}{4\b{[3p/7]}{[p/7]}^2}\mod {p^3} &\t{if $p\e 6\mod
7$.}\endcases$$

\par
 In Section 2, with the help of Maple, we pose many
conjectures on congruences modulo $p^2$ involving the sums
$$\align &\sum_{k=0}^{p-1}\f{\b{2k}k^2}{m^k(k+1)^2},\q
\sum_{k=0}^{p-1}\f{\b{2k}k^2}{m^k(2k-1)^2},
\q\sum_{k=0}^{p-1}\f{\b{2k}k^3}{m^k(k+1)},\q
\sum_{k=0}^{p-1}\f{\b{2k}k^3}{m^k(k+1)^2},\q
\\&\sum_{k=0}^{p-1}\f{\b{2k}k^3}{m^k(2k-1)},\q
\sum_{k=0}^{p-1}\f{\b{2k}k^3}{m^k(2k-1)^2},\q
\sum_{k=0}^{p-1}\f{\b{2k}k^3}{m^k(2k-1)^3},\\&
\sum_{k=0}^{p-1}\f{\b{2k}k^2\b{3k}k}{m^k(k+1)},\q
\sum_{k=0}^{p-1}\f{\b{2k}k^2\b{3k}k}{m^k(k+1)^2},\q
\sum_{k=0}^{p-1}\f{\b{2k}k^2\b{3k}k}{m^k(2k-1)},\q
\\&\sum_{k=0}^{p-1}\f{\b{2k}k^2\b{4k}{2k}}{m^k(k+1)},\q
\sum_{k=0}^{p-1}\f{\b{2k}k^2\b{4k}{2k}}{m^k(k+1)^2}, \q
\sum_{k=0}^{p-1}\f{\b{2k}k^2\b{4k}{2k}}{m^k(2k-1)},
\\&\sum_{k=0}^{p-1}\f{\b{2k}k\b{3k}k\b{6k}{3k}}{m^k(k+1)},\q
\sum_{k=0}^{p-1}\f{\b{2k}k\b{3k}k\b{6k}{3k}}{m^k(k+1)^2},\q
\sum_{k=0}^{p-1}\f{\b{2k}k\b{3k}k\b{6k}{3k}}{m^k(2k-1)}.\endalign$$
As typical examples, we conjecture that for any prime $p>3$,
$$\align
&\sum_{k=0}^{p-1}\f{\b{2k}k^2\b{3k}k}{(-192)^k(k+1)} \e\cases \f
32x^2-4p\mod{p^2}\qq\t{if $3\mid p-1$ and so $4p=x^2+27y^2$,}
\\2(2p+1)\b{[2p/3]}{[p/3]}^2+p\mod {p^2}\q\t{if $p\e 2\mod 3$,}
\endcases
\\&\sum_{k=0}^{p-1}\Ls{\b{2k}k}{2k-1}^3
\e \cases -804y^2-18p\mod {p^2}&\t{if $p=x^2+7y^2\e 1,2,4\mod 7$,}
\\-\f{192}7\b{[3p/7]}{[p/7]}^2\mod p&\t{if $p\e 3\mod 7$,}
\\-\f{3072}{7}\b{[3p/7]}{[p/7]}^2\mod p&\t{if $p\e 5\mod 7$,}
\\-\f{768}{7}\b{[3p/7]}{[p/7]}^2\mod p&\t{if $p\e 6\mod 7$,}
\endcases
\\&\sum_{k=0}^{p-1}\f{\b{2k}k^3}{(-8)^k(2k-1)^2} \e \cases
-4x^2+2p\mod {p^2}&\t{if $p=x^2+4y^2\e 1\mod 4$,}
\\6R_1(p)
\mod {p^2}&\t{if $p\e 3\mod 4$,}\endcases
\\&\Ls p5
\sum_{k=0}^{p-1}\f{\b{2k}k\b{3k}{k}\b{6k}{3k}}{54000^k(k+1)^2}
\\&\q\e\cases \f{2504}{25}x^2-\f{1414+25\sls p5}{25}p\mod {p^2} &\t{if
$p=x^2+3y^2\e 1\mod 3$,}
\\-184R_3(p)+\f{162-25\sls p5}{25}p\mod {p^2}&\t{if $p\e 2\mod 3$ and $p\not=5$,}
\endcases
\\&\sum_{k=0}^{p-1}\f{\b{2k}k^2\b{3k}{k}}{1458^k(2k-1)} \e\cases -\f
{61}{81}\big(4x^2-2p-\f{p^2}{4x^2}\big)-\f
{44}{243}(-1)^{\f{p-1}2}p\mod {p^3}\\\qq\qq\qq\qq \qq\qq\t{if
$p=x^2+3y^2\e 1\mod 3$,}
\\-\f {32}{81}R_3(p)-\f{44}{243}(-1)^{\f{p-1}2}p\mod {p^2}
\q\t{if $p\e 2\mod 3$,}
\endcases
\\&\sum_{k=0}^{p-1}\f{\b{2k}k^2\b{4k}{2k}}{256^k(2k-1)} \e\cases -\f
58(4x^2-2p)\mod {p^2}&\t{if $p=x^2+2y^2\e 1,3\mod 8$,}
\\-\f{5-4(-1)^{\f{p-1}2}}8R_2(p) \mod {p^2}&\t{if $p\e 5,7\mod 8$,}
\endcases
\endalign$$
where
$$\align &R_1(p)=(2p+2-2^{p-1})\b{\f{p-1}2}{\f{p-3}4}^2,
\\&R_3(p)=\Big(1+2p+\f 43(2^{p-1}-1)-\f 32(3^{p-1}-1)\Big)
\b{\f{p-1}2}{[\f p6]}^2,
\\&R_2(p)=\Big(1+(4+2(-1)^{\f{p-1}2})p-4(2^{p-1}-1)-\f p2\sum_{k=1}^{[p/8]}
\f 1k\Big)\b{\f{p-1}2}{[\f p8]}^2.
\endalign$$
\par \par In 2009, Zagier [Z] studied the Ap\'ery-like numbers $\{u_n\}$
satisfying
$$u_0=1,\ u_1=b\qtq{and}(n+1)^2u_{n+1}=(an(n+1)+b)u_n-cn^2u_{n-1}\ (n\ge 1),
$$ where $a,b,c\in\Bbb Z$, $c\not=0$ and $u_n\in\Bbb Z$ for
$n=1,2,3,\ldots$. Let
 $$\align &A'_n=\sum_{k=0}^n\b nk^2\b{n+k}k,\q
 f_n=\sum_{k=0}^n\b nk^3=\sum_{k=0}^n\b nk^2\b{2k}n,\q
 \\&S_n=\sum_{k=0}^{[n/2]}\b{2k}k^2\b n{2k}4^{n-2k}=\sum_{k=0}^n\b nk\b{2k}k\b{2n-2k}{n-k},
 \\&a_n=\sum_{k=0}^n\b nk^2\b{2k}k,
 \q Q_n=\sum_{k=0}^n\b nk(-8)^{n-k}f_k,
 \\&W_n=\sum_{k=0}^{[n/3]}\b{2k}k\b{3k}k\b n{3k}(-3)^{n-3k},
 \\&G_n=\sum_{k=0}^n\b{2k}k^2\b{2n-2k}{n-k}4^{n-k}
 =\sum_{k=0}^n\b nk(-1)^k\b{2k}k^216^{n-k}.\endalign$$
 According to [Z] and [AZ], $\{A'_n\},\ \{f_n\},\ \{S_n\},\ \{a_n\},\
 \{Q_n\}$, $\{W_n\}$ and $\{G_n\}$ are  Ap\'ery-like sequences with
$(a,b,c)=(11,3,-1),(7,2,-8),(12,4,32),(10,3,9),(-17,$ $-6,72), (-9,$
$-3,27)$ and $(32,12,256)$, respectively. The  sequence
 $\{f_n\}$ is called Franel numbers.
 In [S7-S12]  the author
 systematically investigated congruences for sums
 involving $S_n$, $f_n$, $W_n$ and $G_n$.
  For $\{A'_n\},\ \{f_n\},\ \{S_n\},\ \{a_n\},\
 \{Q_n\}$, $\{W_n\}$ and $\{G_n\}$ see A005258, A000172,
 A081085, A002893, A093388,
A291898 and A143583 in Sloane's database ``The On-Line Encyclopedia
of Integer Sequences".
\par Let $p$ be an odd prime.
In [Su3], Z.W. Sun posed many congruences modulo $p^2$ involving
Ap\'ery-like numbers. In [S9,S10,S12], the author conjectured many
congruences modulo $p^3$ involving Ap\'ery-like numbers.
 In Section
3, we pose some conjectures on
$\sum_{k=0}^{p-1}\f{\b{2k}ku_k}{m^k(2k-1)}$ modulo $p^2$, where
$u_n\in\{A'_n,f_n,S_n,a_n,Q_n,W_n,G_n\}$. As typical examples,
$$\align
&\sum_{k=0}^{p-1}\f{\b{2k}kG_k}{64^k(2k-1)} \e
2(-1)^{\f{p-1}2}p^2\mod {p^3},
\\&\sum_{k=0}^{p-1}\f{\b{2k}kS_k}{16^k(2k-1)}
\e\cases 0\mod {p^2}&\t{if $p\e 1\mod 4$,}
\\-R_1(p)\mod {p^2}&\t{if $p\e 3\mod 4$,}
\endcases
\\&
\sum_{k=0}^{p-1}\f{\b{2k}kW_k}{(-12)^k(2k-1)} \e \cases 0\mod
{p^2}&\t{if $p\e 1\mod 3$,}
\\-2(2p+1)\b{[2p/3]}{[p/3]}^2\mod {p^2}&\t{if $p\e 2\mod {3}$}.
\endcases
\\&\sum_{k=0}^{p-1}\f{\b{2k}kf_k}{(-4)^k(2k-1)}
\e\cases 2p-4x^2\mod {p^2}&\t{if $p=x^2+3y^2\e 1\mod 3$,}
\\-8R_3(p)\mod {p^2}&\t{if $p\e 2\mod 3$.}
\endcases
\endalign$$

\section*{2. Conjectures on congruences involving
binomial coefficients}

\par By doing calculations with Maple, we form the following
challenging conjectures.

\pro{Conjecture 2.1} Let $p>3$ be a prime. Then
$$\align &\sum_{k=0}^{p-1}\f{\b{2k}k^2\b{3k}k}{(-192)^k(k+1)}
\e\cases \f 32x^2-4p\mod{p^2}\qq\t{if $3\mid p-1$ and so
$4p=x^2+27y^2$,}
\\2(2p+1)\b{[2p/3]}{[p/3]}^2+p\mod {p^2}\q\t{if $p\e 2\mod 3$,}
\endcases
\\&\sum_{k=0}^{p-1}\f{\b{2k}k^2\b{3k}k}{(-192)^k(k+1)^2}
\e\cases \f {x^2}4-3p\mod{p^2}\qq\t{if $3\mid p-1$ and so
$4p=x^2+27y^2$,}
\\13(2p+1)\b{[2p/3]}{[p/3]}^2+\f p2\mod {p^2}\q\t{if $p\e 2\mod 3$,}
\endcases
\\&\sum_{k=0}^{p-1}\f{\b{2k}k^2\b{3k}{k}}{(-192)^k(2k-1)} \e\cases -\f
34x^2+\f 98p\mod{p^2}\qq\t{if $3\mid p-1$ and $4p=x^2+27y^2$,}
\\-\f 12(2p+1)\b{[2p/3]}{[p/3]}^2+\f 38p\mod {p^2}\q\t{if $p\e 2\mod{3}$.}
\endcases
\endalign$$
\endpro

\pro{Conjecture 2.2} Let $p>3$ be a prime. Then
 $$\align&\Ls
{10}p\sum_{k=0}^{p-1}\f{\b{2k}k\b{3k}k\b{6k}{3k}}{(-12288000)^k(k+1)}
\\&\q\e\cases \f {966}5x^2-\f{3854}5p\mod{p^2}&\t{if $3\mid p-1$ and so
$4p=x^2+27y^2$,}
\\1280(2p+1)\b{[2p/3]}{[p/3]}^2+\f{1922}5p\mod {p^2}&\t{if $p\e 2\mod 3$,}
\endcases
\\&\Ls {10}p\sum_{k=0}^{p-1}\f{\b{2k}k\b{3k}k\b{6k}{3k}}{(-12288000)^k(k+1)^2}
\\&\q\e\cases -\f {112604}{25}x^2+\big(\f{293656}{25}-\ls{10}p\big)p
\mod{p^2}\q\t{if $3\mid p-1$ and so $4p=x^2+27y^2$,}
\\11776(2p+1)\b{[2p/3]}{[p/3]}^2-\big(\f{68448}{25}
+\ls {10}p\big)p\mod {p^2}\q\t{if $p\e 2\mod 3$,}
\endcases
\\&\Ls{10}p\sum_{k=0}^{p-1}\f{\b{2k}k\b{3k}{k}\b{6k}{3k}}
{(-12288000)^k(2k-1)}
\\&\e\cases -\f {177}{200}x^2+\f{53199}{32000}p\mod{p^2}
\qq\t{if $p\e 1\mod 3$ and so $4p=x^2+27y^2$,}
\\-\f 1{20}(2p+1)\b{[2p/3]}{[p/3]}^2+\f{3441}{32000}p
\mod {p^2}\q\t{if $p\e 2\mod{3}$ and $p\not=5$.}
\endcases
\endalign$$
\endpro
\par{\bf Remark 2.1} Let $p$ be a prime with $p>5$. In [S10], the
author conjectured that if $p\e 1\mod 3$ and so $4p=x^2+27y^2$, then
$$\sum_{k=0}^{p-1}\f{\b{2k}k^2\b{3k}k}{(-192)^k}\e
\Ls{10}p\sum_{k=0}^{p-1}\f{\b{2k}k\b{3k}k\b{6k}{3k}}{(-12288000)^k}
\e x^2-2p-\f{p^2}{x^2}\mod {p^3};$$ if $p\e 2\mod 3$, then
$$\sum_{k=0}^{p-1}\f{\b{2k}k^2\b{3k}k}{(-192)^k}
\e\f{800}{161}\Ls{10}p\sum_{k=0}^{p-1}\f{\b{2k}k\b{3k}k\b{6k}{3k}}{(-12288000)^k}
\e\f 34p^2\b{[2p/3]}{[p/3]}^{-2}\mod {p^3}.$$ The congruence for
$\sum_{k=0}^{p-1}\f{\b{2k}k^2\b{3k}k}{(-192)^k}\mod {p^2}$ was
conjectured by Z.W. Sun[Su1] earlier.
\par Let $p>3$ be a prime. In [Su1,Su3], Z.W. Sun conjectured
congruences for $\sum_{k=0}^{p-1}\b{2k}k^3/m^k$ $\mod {p^2}$ with
$m=1,-8,16,-64,256,-512,4096$. Such conjectures were proved by the
author in [S2]. In [S10], the author conjectured congruences for
$\sum_{k=0}^{p-1}\b{2k}k^3/m^k$ $\mod {p^3}$ in the cases
$m=1,-8,16,-64,256,-512,4096$.
 \pro{Conjecture 2.3} Let $p$ be an odd prime. Then
$$\align
&\sum_{k=0}^{p-1}\f{\b{2k}k^3}{(-8)^k(k+1)} \e \f 12R_1(p)+p\mod
{p^2}\qtq{for}p\e 3\mod 4,
\\&\sum_{k=0}^{p-1}\f{\b{2k}k^3}{(-8)^k(k+1)^2}
\e\cases 8x^2-8p\mod {p^2}&\t{if $p=x^2+4y^2\e 1\mod 4$,}
\\3R_1(p)+2p\mod {p^2}&\t{if $p\e 3\mod 4$,}
\endcases
\\&\sum_{k=0}^{p-1}\f{\b{2k}k^3}{(-8)^k(2k-1)}
\e \cases -4x^2\mod {p^2}&\t{if $p=x^2+4y^2\e 1\mod 4$,}
\\2p-2R_1(p)\mod {p^2}&\t{if $p\e 3\mod 4$,}
\endcases
\\&\sum_{k=0}^{p-1}\f{\b{2k}k^3}{(-8)^k(2k-1)^2} \e \cases
-4x^2+2p\mod {p^2}&\t{if $p=x^2+4y^2\e 1\mod 4$,}
\\6R_1(p)
\mod {p^2}&\t{if $p\e 3\mod 4$,}\endcases
\\
&\sum_{k=0}^{p-1}\f{\b{2k}k^3}{(-8)^k(2k-1)^3} \e \cases 48y^2\mod
{p^2}&\t{if $p=x^2+4y^2\e 1\mod 4$,}
\\-6p-12R_1(p)
\mod {p^2}&\t{if $p\e 3\mod 4$.}
\endcases\endalign$$
\endpro

\pro{Conjecture 2.4} Let $p$ be an odd prime. Then
$$\align
&\sum_{k=0}^{p-1}\f{\b{2k}k^3}{64^k(k+1)^2} \e\cases 8x^2-5p\mod
{p^2}&\t{if $p=x^2+4y^2\e 1\mod 4$,}
\\-6R_1(p)-p\mod {p^2}&\t{if $p\e 3\mod 4$,}
\endcases
\\&\sum_{k=0}^{p-1}\f{\b{2k}k^3}{64^k(2k-1)} \e \cases p-2x^2\mod
{p^2}&\t{if $p=x^2+4y^2\e 1\mod 4$,}
\\-\f 12R_1(p)\mod {p^2}&\t{if $p\e 3\mod 4$,}
\endcases
\\&\sum_{k=0}^{p-1}\f{\b{2k}k^3}{64^k(2k-1)^2}
\e \cases 2x^2-p\mod {p^2}&\t{if $p=x^2+4y^2\e 1\mod 4$,}
\\\f 32R_1(p)\mod {p^2}&\t{if $p\e 3\mod 4$,}
\endcases
\\&\sum_{k=0}^{p-1}\f{\b{2k}k^3}{64^k(2k-1)^3}
\e \cases 0\mod {p^2}&\t{if $p=x^2+4y^2\e 1\mod 4$,}
\\-3R_1(p)\mod {p^2}&\t{if $p\e 3\mod 4$.}
\endcases
\endalign $$
\endpro
\pro{Conjecture 2.5} Let $p$ be an odd prime. Then
$$\align
&(-1)^{[\f p4]}\sum_{k=0}^{p-1}\f{\b{2k}k^3}{(-512)^k(k+1)} \e\cases
8x^2-6p\mod {p^2}&\t{if $p=x^2+4y^2\e 1\mod 4$,}
\\-4R_1(p)-2p\mod {p^2}&\t{if $p\e 3\mod 4$,}
\endcases
\\&(-1)^{[\f p4]}\sum_{k=0}^{p-1}\f{\b{2k}k^3}{(-512)^k(k+1)^2}
\\&\qq\e\cases -16x^2+(7-(-1)^{[\f p4]})p\mod {p^2}&\t{if $p=x^2+4y^2\e
1\mod 4$,}
\\-24R_1(p)-(-1)^{[\f p4]}p\mod {p^2}&\t{if $p\e 3\mod 4$.}
\endcases
\\&(-1)^{[\f p4]}\sum_{k=0}^{p-1}\f{\b{2k}k^3}{(-512)^k(2k-1)}
\e \cases \f 54p-3x^2\mod {p^2}&\t{if $p=x^2+4y^2\e 1\mod 4$,}
\\-\f p4+\f{1}4R_1(p)\mod {p^2}
&\t{if $4\mid p-3$,}
\endcases
\\&(-1)^{[\f p4]}\sum_{k=0}^{p-1}\f{\b{2k}k^3}{(-512)^k(2k-1)^2}
\e \cases 2x^2-\f 58p\mod {p^2}\qq\t{if $p=x^2+4y^2\e 1\mod 4$,}
\\-\f 34R_1(p)+\f 38p\mod {p^2}
\q\t{if $p\e 3\mod 4$,}
\endcases
\\&(-1)^{[\f p4]}\sum_{k=0}^{p-1}\f{\b{2k}k^3}{(-512)^k(2k-1)^3}
\e \cases -\f 32x^2+\f 38p\mod {p^2}&\t{if $p=x^2+4y^2\e 1\mod 4$,}
\\\f 32R_1(p)-\f 38p\mod {p^2}
&\t{if $4\mid p-3$.}
\endcases
\endalign$$
\endpro

\pro{Conjecture 2.6} Let $p>3$ be a prime. Then
$$\align
&\sum_{k=0}^{p-1}\f{\b{2k}k^2\b{4k}{2k}}{648^k(k+1)} \e\cases \f
{10}3x^2-\f 43p\mod {p^2}&\t{if $p=x^2+4y^2\e 1\mod 4$,}
\\-\f 32R_1(p)-\f p3\mod {p^2}&\t{if $p\e 3\mod 4$,}
\endcases
\\&\sum_{k=0}^{p-1}\f{\b{2k}k^2\b{4k}{2k}}{648^k(k+1)^2}
\e\cases \f {112}9x^2-\f {64}9p\mod {p^2}&\t{if $p=x^2+4y^2\e 1\mod
4$,}
\\-11R_1(p)-\f{10}9 p\mod {p^2}&\t{if $p\e 3\mod 4$,}
\endcases
\\&\sum_{k=0}^{p-1}\f{\b{2k}k^2\b{4k}{2k}}{648^k(2k-1)}
\e \cases -\f{76}{27}x^2+\f{104}{81}p\mod {p^2}&\t{if $p=x^2+4y^2\e
1\mod 4$,}
\\-\f 29R_1(p)+\f{10}{81}p
\mod {p^2}&\t{if $p\e 3\mod 4$.}
\endcases\endalign$$
\endpro

\pro{Conjecture 2.7} Let $p>3$ be a prime. Then
$$\align &\Ls p3\sum_{k=0}^{p-1}\f{\b{2k}k\b{3k}k\b{6k}{3k}}{12^{3k}(k+1)}
\e\f 35R_1(p)\mod {p^2}\qtq{for}p\e 3\mod 4,
\\&\Ls p3\sum_{k=0}^{p-1}\f{\b{2k}k\b{3k}k\b{6k}{3k}}{12^{3k}(k+1)^2}
\e\cases 8x^2-(4+\sls p3)p\mod {p^2}&\t{if $p=x^2+4y^2\e 1\mod 4$,}
\\\f{138}{25}R_1(p)-\sls p3p\mod {p^2}&\t{if $p\e 3\mod 4$,}
\endcases
\\&\Ls
p3\sum_{k=0}^{p-1}\f{\b{2k}k\b{3k}{k}\b{6k}{3k}}{12^{3k}(2k-1)} \e
\cases -\f{26}{9}x^2+\f{13}{9}p\mod {p^2}&\t{if $p=x^2+4y^2\e 1\mod
4$,}
\\\f 16R_1(p)\mod {p^2}
&\t{if $p\e 3\mod 4$.}
\endcases\endalign$$
\endpro

\pro{Conjecture 2.8} Let $p$ be a prime with $p\not=2,3,11$. Then
$$\align
&\Ls {33}p\sum_{k=0}^{p-1}\f{\b{2k}k\b{3k}k\b{6k}{3k}}{66^{3k}(k+1)}
\e\cases -26x^2+28p\mod {p^2}&\t{if $p=x^2+4y^2\e 1\mod 4$,}
\\-\f {363}{10}R_1(p)-15p\mod {p^2}&\t{if $p\e 3\mod 4$,}
\endcases
\\&\Ls {33}p\sum_{k=0}^{p-1}\f{\b{2k}k\b{3k}k\b{6k}{3k}}{66^{3k}(k+1)^2}
\\&\q\e\cases 488x^2-\big(295+\sls {33}p\big)p\mod {p^2}&\t{if
$p=x^2+4y^2\e 1\mod 4$,}
\\-\f{8349}{25}R_1(p)+\big(51-\sls {33}p\big)p\mod {p^2}
&\t{if $p\e 3\mod 4$,}
\endcases
\\&\Ls {33}p\sum_{k=0}^{p-1}\f{\b{2k}k\b{3k}{k}\b{6k}{3k}}{66^{3k}(2k-1)}
\e \cases -\f{3716}{1089}x^2+\f{18848}{11979}p\mod {p^2}&\t{if
$p=x^2+4y^2\e 1\mod 4$,}
\\-\f 2{33}R_1(p)
+\f{530}{3993}p\mod {p^2}&\t{if $p\e 3\mod 4$}.
\endcases\endalign$$
\endpro

\pro{Conjecture 2.9} Let $p>3$ be a prime. Then
$$\align &\sum_{k=0}^{p-1}\f{\b{2k}k^3}{16^k(k+1)}
\e\cases \f{16}3x^2-\f{10}3p\mod {p^2}&\t{if $p=x^2+3y^2\e 1\mod
3$,}
\\-\f 43R_3(p)-\f 23p\mod {p^2}&\t{if $p\e 2\mod 3$,}
\endcases
\\&\sum_{k=0}^{p-1}\f{\b{2k}k^3}{16^k(k+1)^2}
\e\cases 8x^2-7p\mod {p^2}&\t{if $p=x^2+3y^2\e 1\mod 3$,}
\\-8R_3(p)-3p\mod {p^2}&\t{if $p\e 2\mod 3$,}
\endcases\\&\sum_{k=0}^{p-1}\f{\b{2k}k^3}{16^k(2k-1)}
\e \cases 4y^2\mod {p^2}&\t{if $p=x^2+3y^2\e 1\mod 3$,}
\\-\f{8}3R_3(p)+\f 23p\mod {p^2}&\t{if $p\e 2\mod 3$,}
\endcases
\\ &\sum_{k=0}^{p-1}\f{\b{2k}k^3}{16^k(2k-1)^2}\e
\cases 4x^2-2p\mod {p^2}&\t{if $p=x^2+3y^2\e 1\mod 3$,}
\\ 8R_3(p)\mod {p^2}&\t{if $p\e 2\mod
3$,}\endcases
\\&\sum_{k=0}^{p-1}\f{\b{2k}k^3}{16^k(2k-1)^3}\e
\cases -12y^2 \mod {p^2} &\t{if $p=x^2+3y^2\e 1\mod 3$,}
\\-16R_3(p)-2p\mod{p^2}
&\t{if $p\e 2\mod 3$.}
\endcases\endalign$$
\endpro
\pro{Conjecture 2.10} Let $p>3$ be a prime. Then
$$\align
&(-1)^{\f{p-1}2}\sum_{k=0}^{p-1}\f{\b{2k}k^3}{256^k(k+1)} \e\cases
\f{8}3x^2-\f{2}3p\mod {p^2}&\t{if $p=x^2+3y^2\e 1\mod 3$,}
\\\f {16}3R_3(p)+\f 23p\mod {p^2}&\t{if $p\e 2\mod 3$,}
\endcases
\\&(-1)^{\f{p-1}2}\sum_{k=0}^{p-1}\f{\b{2k}k^3}{256^k(k+1)^2}
\\&\q\e\cases 16x^2-(8+(-1)^{\f{p-1}2})p\mod {p^2}&\t{if $p=x^2+3y^2\e
1\mod 3$,}
\\32R_3(p)-(-1)^{\f{p-1}2}p\mod {p^2}&\t{if $p\e 2\mod 3$,}
\endcases
\\&(-1)^{\f{p-1}2}\sum_{k=0}^{p-1}\f{\b{2k}k^3}{256^k(2k-1)}
\e \cases 8y^2-\f 32p\mod {p^2}&\t{if $p=x^2+3y^2\e 1\mod 3$,}
\\\f{2}3R_3(p)-\f p6\mod {p^2}&\t{if $p\e 2\mod 3$,}
\endcases
\\&(-1)^{\f{p-1}2}\sum_{k=0}^{p-1}\f{\b{2k}k^3}{256^k(2k-1)^2}
\e\cases 2x^2-\f 34p\mod {p^2}&\t{if $p=x^2+3y^2\e 1\mod 3$,}
\\-2R_3(p)+\f p4\mod {p^2}&\t{if $p\e 2\mod 3$,}
\endcases
\\&(-1)^{\f{p-1}2}\sum_{k=0}^{p-1}\f{\b{2k}k^3}{256^k(2k-1)^3}
\e\cases \f p4-x^2 \mod {p^2}&\t{if $p=x^2+3y^2\e 1\mod 3$,}
\\4R_3(p)-\f p4\mod {p^2}&\t{if $p\e 2\mod 3$.}
\endcases\endalign$$
Moreover,
 $$(-1)^{\f{p-1}2}\sum_{k=0}^{p-1}\f{\b{2k}k^3}{256^k(2k-1)} \e
-\f 14\sum_{k=0}^{p-1}\f{\b{2k}k^3}{16^k(2k-1)}\mod
{p^3}\qtq{for}p\e 2\mod 3.$$
 \endpro

\pro{Conjecture 2.11} Let $p>3$ be a prime. Then
$$\align &\sum_{k=0}^{p-1}\f{\b{2k}k^2\b{3k}{k}}{108^k(k+1)}
\e -2R_3(p)\mod {p^2}\qtq{for}p\e 2\mod 3,
\\&\sum_{k=0}^{p-1}\f{\b{2k}k^2\b{3k}{k}}{108^k(k+1)^2}
\e\cases 8x^2-5p\mod {p^2}&\t{if $p=x^2+3y^2\e 1\mod 3$,}
\\-13R_3(p)-p\mod {p^2}&\t{if $p\e 2\mod 3$,}
\endcases\\&\sum_{k=0}^{p-1}\f{\b{2k}k^2\b{3k}{k}}{108^k(2k-1)}
\e\cases -\f 59(4x^2-2p)\mod {p^2}&\t{if $p=x^2+3y^2\e 1\mod 3$,}
\\-\f 89R_3(p)\mod {p^2}&\t{if $p\e 2\mod 3$,}
\endcases\endalign$$
\endpro

\pro{Conjecture 2.12} Let $p>3$ be a prime. Then
$$\align&\sum_{k=0}^{p-1}\f{\b{2k}k^2\b{4k}{2k}}{(-144)^k(k+1)} \e\cases
\f{16}3x^2-\f{10}3p\mod {p^2}&\t{if $p=x^2+3y^2\e 1\mod 3$,}
\\\f 43R_3(p)+\f 23p\mod {p^2}&\t{if $p\e 2\mod 3$,}
\endcases
\\&\sum_{k=0}^{p-1}\f{\b{2k}k^2\b{4k}{2k}}{(-144)^k(k+1)^2}
\e\cases \f{40}9x^2-\f{43}9p\mod {p^2}&\t{if $p=x^2+3y^2\e 1\mod
3$,}
\\\f {88}9R_3(p)+\f 59p\mod {p^2}&\t{if $p\e 2\mod 3$,}
\endcases
\\&\sum_{k=0}^{p-1}\f{\b{2k}k^2\b{4k}{2k}}{(-144)^k(2k-1)} \e\cases
-\f {28}9x^2+\f 89p\mod {p^2}&\t{if $p=x^2+3y^2\e 1\mod 3$,}
\\-\f 89R_3(p)+\f 23p\mod {p^2}&\t{if $p\e 2\mod 3$}
\endcases
\endalign$$
\endpro

\pro{Conjecture 2.13} Let $p>5$ be a prime. Then
$$\align &\Ls p5
\sum_{k=0}^{p-1}\f{\b{2k}k\b{3k}{k}\b{6k}{3k}}{54000^k(k+1)}
\e\cases -\f{16}5x^2+\f{26}5p\mod {p^2} &\t{if $p=x^2+3y^2\e 1\mod
3$,}
\\-20R_3(p)-\f{18}5p\mod {p^2}&\t{if $p\e 2\mod 3$,}
\endcases
\\&\Ls p5
\sum_{k=0}^{p-1}\f{\b{2k}k\b{3k}{k}\b{6k}{3k}}{54000^k(k+1)^2}
\\&\q\e\cases \f{2504}{25}x^2-\f{1414+25\sls p5}{25}p\mod {p^2} &\t{if
$p=x^2+3y^2\e 1\mod 3$,}
\\-184R_3(p)+\f{162-25\sls p5}{25}p\mod {p^2}&\t{if $p\e 2\mod 3$,}
\endcases\\&\Ls p5\sum_{k=0}^{p-1}\f{\b{2k}k\b{3k}{k}\b{6k}{3k}}{54000^k(2k-1)}
\e\cases -\f {748}{225}x^2+\f{1708}{1125}p\mod {p^2}&\t{if
$p=x^2+3y^2\e 1\mod 3$,}
\\-\f 8{45}R_3(p)+\f{18}{125}p\mod {p^2}&\t{if $p\e 2\mod 3$.}
\endcases\endalign$$
\endpro

\pro{Conjecture 2.14} Let $p>3$ be a prime. Then
$$\align &
\sum_{k=0}^{p-1}\f{\b{2k}k^2\b{3k}{k}}{1458^k(k+1)} \e\cases
2(-1)^{\f{p-1}2}p-p^2\mod {p^3}&\t{if $p\e 1\mod 3$,}
\\-12R_3(p)+2(-1)^{\f{p-1}2}p\mod {p^2}&\t{if $p\e 2\mod 3$,}
\endcases
\\&\sum_{k=0}^{p-1}\f{\b{2k}k^2\b{3k}{k}}{1458^k(k+1)^2}
\\&\q\e\cases 42x^2-(22+2(-1)^{\f{p-1}2})p\mod {p^2}&\t{if $p=x^2+3y^2\e
1\mod 3$,}
\\-78R_3(p)-(1+2(-1)^{\f{p-1}2})p\mod {p^2}&\t{if $p\e 2\mod 3$,}
\endcases
\\&\sum_{k=0}^{p-1}\f{\b{2k}k^2\b{3k}{k}}{1458^k(2k-1)}
\e\cases -\f {61}{81}\big(4x^2-2p-\f{p^2}{4x^2}\big)-\f
{44}{243}(-1)^{\f{p-1}2}p\mod {p^3}\\\qq\qq\qq\ \qq\qq\t{if
$p=x^2+3y^2\e 1\mod 3$,}
\\-\f {32}{81}R_3(p)-\f{44}{243}(-1)^{\f{p-1}2}p\mod {p^2}
\q\t{if $p\e 2\mod 3$.}
\endcases\endalign$$
\endpro

\pro{Conjecture 2.15} Let $p$ be an odd prime. Then
$$\align
&(-1)^{\f{p-1}2}\sum_{k=0}^{p-1}\f{\b{2k}k^3}{(-64)^k(k+1)} \e\cases
6x^2-4p\mod {p^2}\qq\t{if $p=x^2+2y^2\e 1,3\mod 8$,}
\\-\f {5-4(-1)^{\f{p-1}2}}2R_2(p)-p\mod {p^2}\q\t{if $p\e 5,7\mod 8$,}
\endcases
\\&(-1)^{\f{p-1}2}\sum_{k=0}^{p-1}\f{\b{2k}k^3}{(-64)^k(k+1)^2}
\cases 4x^2-5p\mod {p^2}&\t{if $p=x^2+2y^2\e 1\mod 8$,}
\\4x^2-3p\mod {p^2}&\t{if $p=x^2+2y^2\e 3\mod 8$,}
\\-3R_2(p)-3p\mod {p^2}&\t{if $p\e 5\mod 8$,}
\\-27R_2(p)-p\mod {p^2}&\t{if $p\e 7\mod 8$,}
\endcases
\\&(-1)^{\f{p-1}2}\sum_{k=0}^{p-1}\f{\b{2k}k^3}{(-64)^k(2k-1)}
\e \cases p-3x^2\mod {p^2}\qq\t{if $p=x^2+2y^2\e 1,3\mod 8$,}
\\\f{5-4(-1)^{\f{p-1}2}}4R_2(p)-\f p2\mod {p^2}\q\t{if $p\e 5,7\mod 8$,}
\endcases
\\ &(-1)^{\f{p-1}2}\sum_{k=0}^{p-1}\f{\b{2k}k^3}{(-64)^k(2k-1)^2}
\e\cases x^2\mod{p^2}\qq\qq\q\t{if $p=x^2+2y^2\e 1,3\mod 8$,}
\\-\f 34(5-4(-1)^{\f{p-1}2})R_2(p)+\f p2
\mod{p^2}\ \t{if $p\e 5,7\mod 8$,}
\endcases
\\&(-1)^{\f{p-1}2}\sum_{k=0}^{p-1}\f{\b{2k}k^3}{(-64)^k(2k-1)^3}
\e\cases p-2x^2\mod{p^2}\qq\t{if $p=x^2+2y^2\e 1,3\mod 8$,}
\\ \f 32(5-4(-1)^{\f{p-1}2})R_2(p)\mod{p^2}\q\t{if $p\e 5,7\mod 8$.}
\endcases
\endalign$$
\endpro

\pro{Conjecture 2.16} Let $p$ be an odd prime. Then
$$\align &\sum_{k=0}^{p-1}\f{\b{2k}k^2\b{4k}{2k}}{256^k(k+1)}
\e \cases -\f 13R_2(p)\mod {p^2}&\t{if $p\e 5\mod 8$,}
\\-3R_2(p)\mod {p^2}&\t{if $p\e  7\mod 8$,}
\endcases
\\&\sum_{k=0}^{p-1}\f{\b{2k}k^2\b{4k}{2k}}{256^k(k+1)^2} \e \cases
8x^2-5p\mod {p^2}&\t{if $p=x^2+2y^2\e 1,3\mod 8$,}
\\ -\f
{22}9R_2(p)-p\mod {p^2}&\t{if $p\e  5\mod 8$,}
\\-22R_2(p)-p\mod {p^2}&\t{if $p\e  7\mod 8$,}
\endcases \\&\sum_{k=0}^{p-1}\f{\b{2k}k^2\b{4k}{2k}}{256^k(2k-1)} \e\cases -\f
58(4x^2-2p)\mod {p^2}&\t{if $p=x^2+2y^2\e 1,3\mod 8$,}
\\-\f{5-4(-1)^{\f{p-1}2}}8R_2(p) \mod {p^2}&\t{if $p\e 5,7\mod 8$}
\endcases\endalign$$
\endpro

\pro{Conjecture 2.17} Let $p$ be a prime with $p\not=2,7$. Then
 $$\align
&\sum_{k=0}^{p-1}\f{\b{2k}k^2\b{4k}{2k}}{28^{4k}(k+1)}
\\&\q\e\cases -284x^2+\big(142+144\ls p3\big)p\mod {p^2}&\t{if
$p=x^2+2y^2\e 1,3\mod 8$,}
\\-147R_2(p)+144\ls p3p\mod {p^2}&\t{if $p\e 5\mod 8$,}
\\-1323R_2(p)+144\ls p3p\mod {p^2}&\t{if $p\e 7\mod 8$,}
\endcases
\\&\sum_{k=0}^{p-1}\f{\b{2k}k^2\b{4k}{2k}}{28^{4k}(k+1)^2}
\\&\q\e\cases 5576x^2-\big(3652+\ls p3\big)p\mod {p^2}&\t{if
$p=x^2+2y^2\e 1\mod 8$,}
\\5576x^2-\big(2789+864\ls p3\big)p\mod {p^2}&\t{if
$p=x^2+2y^2\e 3\mod 8$,}
\\-1078R_2(p)-\big(1+864\ls p3\big)p\mod {p^2}&\t{if $p\e 5\mod 8$,}
\\-9702R_2(p)-\big(1+864\ls p3\big)p\mod {p^2}&\t{if $p\e 7\mod 8$,}
\endcases
 \\& \sum_{k=0}^{p-1}\f{\b{2k}k^2\b{4k}{2k}}{28^{4k}(2k-1)}
\\&\e\cases -\f{2363}{686}x^2+\f{16541-1224\sls p3}{9604}p\mod
{p^2}&\t{if $p=x^2+2y^2\e 1,3\mod 8$,}
\\-\f {45-36(-1)^{(p-1)/2}}{392}R_2(p)-\f{306}{2401}\Ls p3p
\mod {p^2}&\t{if $p\e 5,7\mod 8$.}
\endcases\endalign$$
\endpro

\pro{Conjecture 2.18} Let $p$ be an odd prime. Then
$$\align
&\sum_{k=0}^{p-1}\f{\b{2k}k^2\b{3k}k}{8^k(k+1)} \e\cases
\f{11}2x^2-\f 72p\mod {p^2}\qq\t{if $p=x^2+2y^2\e 1,3\mod 8$,}
\\-\f{5-4(-1)^{\f{p-1}2}}8R_2(p)-\f 34p\mod {p^2}\q\t{if $p\e 5,7\mod
8$,}
\endcases
\\&\sum_{k=0}^{p-1}\f{\b{2k}k^2\b{3k}k}{8^k(k+1)^2}
\e\cases \f{31}4x^2-\f {29}4p\mod {p^2}&\t{if $p=x^2+2y^2\e 1,3\mod
8$,}
\\-\f{13}{16}R_2(p)-\f {27}8p\mod {p^2}&\t{if $p\e 5\mod
8$,}
\\-\f{117}{16}R_2(p)-\f {27}8p\mod {p^2}&\t{if $p\e 7\mod
8$,}
\endcases
\\&\sum_{k=0}^{p-1}\f{\b{2k}k^2\b{3k}{k}}{8^k(2k-1)} \e \cases
\f{7}{2}p-x^2\mod {p^2}\qq\qq\t{if $p=x^2+2y^2\e 1,3\mod 8$,}
\\-\f 34(5-4(-1)^{\f{p-1}2})R_2(p)+3p\mod {p^2}\q\t{if $p\e 5,7\mod 8$,}
\endcases\endalign$$\endpro

\pro{Conjecture 2.19} Let $p$ be a prime with $p\not=2,5$. Then
$$\align
&\Ls{-5}p\sum_{k=0}^{p-1}\f{\b{2k}k\b{3k}k\b{6k}{3k}}{20^{3k}(k+1)}
 \e\cases
\f{14}5x^2-\f 45p\mod {p^2}&\t{if $p=x^2+2y^2\e 1,3\mod 8$,}
\\\f 56R_2(p)+\f 35p\mod {p^2}&\t{if $p\e 5\mod
8$,}\\\f{15}2R_2(p)+\f 35p\mod {p^2}&\t{if $p\e 7\mod 8$,}
\endcases
\\&\Ls{-5}p\sum_{k=0}^{p-1}\f{\b{2k}k\b{3k}k\b{6k}{3k}}{20^{3k}(k+1)^2}
\\&\q\e\cases \f{484}{25}x^2-\big(\f{244}{25}+\ls{-5}p\big)p\mod
{p^2}&\t{if $p=x^2+2y^2\e 1,3\mod 8$,}
\\\f{23}3R_2(p)-\big(\f{2}{25}+\ls{-5}p\big)p\mod {p^2}&\t{if $p\e 5\mod
8$,}
\\69R_2(p)-\big(\f{2}{25}+\ls{-5}p\big)p\mod {p^2}&\t{if $p\e 7\mod
8$,}
\endcases
\\&\Ls{-5}p\sum_{k=0}^{p-1}\f{\b{2k}k\b{3k}k\b{6k}{3k}}{20^{3k}(2k-1)}
\e\cases -\f {79}{25}x^2+\f{181}{125}p\mod {p^2}\qq\t{if
$p=x^2+2y^2\e 1,3\mod 8$,}
\\\f{5-4(-1)^{\f{p-1}2}}{20}R_2(p) -\f{33}{250}p \mod {p^2}
\q\t{if $p\e 5,7\mod 8$.}
\endcases\endalign$$
\endpro
\par {\bf Remark 2.2}  Let $p$ be a prime with $p>7$ and $p\not=71$.
In [S10], the author conjectured that if $p\e 1,3\mod 8$ and so
$p=x^2+2y^2$, then
$$\align \sum_{k=0}^{p-1}\f{\b{2k}k^2\b{4k}{2k}}{256^k}
&\e \sum_{k=0}^{p-1}\f{\b{2k}k^2\b{4k}{2k}}{28^{4k}} \e \Ls
{-5}p\sum_{k=0}^{p-1}\f{\b{2k}k\b{3k}k\b{6k}{3k}}{20^{3k}}\e
4x^2-2p-\f {p^2}{4x^2}\mod{p^3};\endalign$$ if $p\e 5,7\mod 8$, then
$$\align \sum_{k=0}^{p-1}\f{\b{2k}k^2\b{4k}{2k}}{256^k}
&\e -\f{441}{71}\sum_{k=0}^{p-1}\f{\b{2k}k^2\b{4k}{2k}}{28^{4k}} \e
-\f{25}7\Ls
{-5}p\sum_{k=0}^{p-1}\f{\b{2k}k\b{3k}k\b{6k}{3k}}{20^{3k}}
\\&\e\cases\f{p^2}3\b{[p/4]}{[p/8]}^{-2}\mod
{p^3}&\t{if $p\e 5\mod 8$,}
\\-\f 32p^2\b{[p/4]}{[p/8]}^{-2}\mod {p^3}&\t{if $p\e 7\mod 8$.}
\endcases\endalign$$
The conjecture for $\sum_{k=0}^{p-1}\f{\b{2k}k^2\b{4k}{2k}}{28^{4k}}
\mod {p^2}$ was due to Z.W. Sun[Su1].
\par For any odd prime $p$ let
$$R_7(p)=\sum_{k=0}^{p-1}\f{\b{2k}k^3}{k+1}.$$
\pro{Conjecture 2.20} Let $p$ be a prime with $p\not=2,7$. Then
$$\align &R_7(p)=\sum_{k=0}^{p-1}\f{\b{2k}k^3}{k+1}
\e \cases -44y^2+2p\mod {p^2}&\t{if $p=x^2+7y^2\e 1,2,4\mod 7$,}
\\-\f{1}7\b{[3p/7]}{[p/7]}^2\mod p&\t{if $p\e 3\mod 7$,}
\\-\f{16}7\b{[3p/7]}{[p/7]}^2\mod p&\t{if $p\e 5\mod 7$,}
\\-\f{4}7\b{[3p/7]}{[p/7]}^2\mod p&\t{if $p\e 6\mod 7$,}
\endcases
\\&\sum_{k=0}^{p-1}\f{\b{2k}k^3}{(k+1)^2}
\e \cases -68y^2\mod {p^2}&\t{if $p=x^2+7y^2\e 1,2,4\mod 7$,}
\\6R_7(p)+2p\mod {p^2}&\t{if $p\e 3,5,6\mod 7$,}
\endcases
\\&\sum_{k=0}^{p-1}\f{\b{2k}k^3}{2k-1} \e \cases 14p-36y^2\mod
{p^2}\q\t{if $p=x^2+7y^2\e 1,2,4\mod 7$,}
\\32R_7(p)+48p\mod {p^2}\q\t{if $p\e 3,5,6\mod 7$,}
\endcases
\\ &\sum_{k=0}^{p-1}\f{\b{2k}k^3}{(2k-1)^2}\e \cases -284y^2+34p\mod
{p^2}\q\t{if $p=x^2+7y^2\e 1,2,4\mod 7$,}
\\ -96R_7(p)-96p\mod {p^2}\q\t{if $p\e
3,5,6\mod 7$,}\endcases
\\& \sum_{k=0}^{p-1}\f{\b{2k}k^3}{(2k-1)^3}\e
\cases -804y^2-18p\mod {p^2}&\t{if $p=x^2+7y^2\e 1,2,4\mod 7$,} \\
192R_7(p)+144p\mod {p^2}&\t{if $p\e 3,5,6\mod 7$.}\endcases
\endalign$$
\endpro

\pro{Conjecture 2.21} Let $p$ be a prime with $p\not=2,7$. Then
$$\align
&(-1)^{\f{p-1}2}\sum_{k=0}^{p-1}\f{\b{2k}k^3}{4096^k(k+1)} \\&\q\e
\cases 72y^2+2p\mod {p^2}&\t{if $p=x^2+7y^2\e 1,2,4\mod 7$,}
\\-64R_7(p)-66p\mod {p^2}&\t{if $p\e 3,5,6\mod 7$,}
\endcases
\\&(-1)^{\f{p-1}2}\sum_{k=0}^{p-1}\f{\b{2k}k^3}{4096^k(k+1)^2}
\\&\q\e \cases
-1136y^2+(64-(-1)^{\f{p-1}2})p\mod {p^2}&\t{if $p=x^2+7y^2\e
1,2,4\mod 7$,}
\\-384R_7(p)-(456+(-1)^{\f{p-1}2})p\mod {p^2}&\t{if $p\e 3,5,6\mod 7$,}
\endcases
\\&(-1)^{\f{p-1}2}\sum_{k=0}^{p-1}\f{\b{2k}k^3}{4096^k(2k-1)}
\\&\q\e \cases -\f 74p+22y^2\mod {p^2}&\t{if $p=x^2+7y^2\e 1,2,4\mod 7$,}
\\-\f 12R_7(p)-\f 34p\mod {p^2}&\t{if $p\e 3,5,6\mod 7$,}
\endcases
\\&(-1)^{\f{p-1}2}\sum_{k=0}^{p-1}\f{\b{2k}k^3}{4096^k(2k-1)^2}
\\&\q\e\cases -17y^2+\f{97}{64}p\mod{p^2}&\t{if $p=x^2+7y^2\e 1,2,4\mod
7$,}
\\\f 32R_7(p)-\f{3}{64}p\mod{p^2}&\t{if $p\e 3,5,6\mod 7$,}
\endcases
\\& (-1)^{\f{p-1}2}\sum_{k=0}^{p-1}\f{\b{2k}k^3}{4096^k(2k-1)^3}
\\&\q\e\cases \f{201}{16}y^2+\f{81}{64}\mod{p^2} \q\t{if $p=x^2+7y^2\e
1,2,4\mod 7$,}
\\-3R_7(p)-\f{243}{64}p\mod{p^2}\q\t{if $p\e 3,5,6\mod 7$}
\endcases
\endalign$$
and
$$(-1)^{\f{p-1}2}\sum_{k=0}^{p-1}\f{\b{2k}k^3}{4096^k(2k-1)}
\e -\f 1{64}\sum_{k=0}^{p-1}\f{\b{2k}k^3}{2k-1}\mod{p^3}\qtq{for}
p\e 3,5,6\mod 7.$$
\endpro

\pro{Conjecture 2.22} Let $p$ be a prime with $p\not=2,3,7$. Then
$$\align
&\sum_{k=0}^{p-1}\f{\b{2k}k^2\b{4k}{2k}}{81^k(k+1)} \e \cases
-\f{100}{3}y^2+2p\mod {p^2}&\t{if $p=x^2+7y^2\e 1,2,4\mod 7$,}
\\\f{3}2R_7(p)+\f 43p\mod {p^2}&\t{if $p\e 3,5,6\mod 7$,}
\endcases
\\&\sum_{k=0}^{p-1}\f{\b{2k}k^2\b{4k}{2k}}{81^k(k+1)^2}
\e \cases -\f{436}{9}y^2+\f 43p\mod {p^2}&\t{if $p=x^2+7y^2\e
1,2,4\mod 7$,}
\\11R_7(p)+\f{94}{9}p\mod {p^2}&\t{if $p\e 3,5,6\mod 7$,}
\endcases
\\&\sum_{k=0}^{p-1}\f{\b{2k}k^2\b{4k}{2k}}{81^k(2k-1)}
\e \cases \f{436}{27}y^2-\f{50}{81}p\mod {p^2}&\t{if $p=x^2+7y^2\e
1,2,4\mod 7$,}
\\\f{16}9R_7(p)+\f{208}{81}p\mod {p^2}&\t{if $p\e 3,5,6\mod 7$.}
\endcases\endalign$$
\endpro

\pro{Conjecture 2.23} Let $p$ be a prime with $p\not=2,3,7$. Then
$$\align
& \sum_{k=0}^{p-1}\f{\b{2k}k^2\b{4k}{2k}} {(-3969)^k(k+1)} \e \cases
\f{28}3x^2-\f{22}{3}p\mod {p^2}&\t{if $p=x^2+7y^2\e 1,2,4\mod 7$,}
\\-21R_7(p)-\f{64}{3}p\mod {p^2}&\t{if $p\e 3,5,6\mod 7$,}
\endcases
\\& \sum_{k=0}^{p-1}\f{\b{2k}k^2\b{4k}{2k}} {(-3969)^k(k+1)^2} \e \cases
-\f{356}9x^2+\f{188}{9}p\mod {p^2}&\t{if $p=x^2+7y^2\e 1,2,4\mod
7$,}
\\-154R_7(p)-\f{1612}{9}p\mod {p^2}&\t{if $p\e 3,5,6\mod 7$,}
\endcases
 \\& \sum_{k=0}^{p-1}\f{\b{2k}k^2\b{4k}{2k}}
{(-3969)^k(2k-1)} \e \cases \f{2102}{189}y^2-\f{7090}{3969}p\mod
{p^2}&\t{if $p=x^2+7y^2\e 1,2,4\mod 7$,}
\\\f{32}{63}R_7(p)+\f{3088}{3969}p\mod {p^2}&\t{if $p\e 3,5,6\mod 7$.}
\endcases\endalign$$

\endpro

\pro{Conjecture 2.24} Let $p$ be a prime with $p>7$. Then
$$\align&\Ls {-15}p\sum_{k=0}^{p-1}\f{\b{2k}k\b{3k}k\b{6k}{3k}}
{(-15)^{3k}(k+1)}\\&\q \e \cases -\f{188}{5}y^2+2p\mod {p^2}&\t{if
$p=x^2+7y^2\e 1,2,4\mod 7$,}
\\\f {15}{128}R_7(p)-\f{81}{40}p\mod {p^2}&\t{if $p\e 3,5,6\mod 7$,}
\endcases
\\&\Ls {-15}p\sum_{k=0}^{p-1}\f{\b{2k}k\b{3k}k\b{6k}{3k}}
{(-15)^{3k}(k+1)^2}\\&\q \e \cases \f{172}{25}y^2-\big(\f
75+\ls{-15}p\big)p\mod {p^2}&\t{if $p=x^2+7y^2\e 1,2,4\mod 7$,}
\\\f {69}{64}R_7(p)-\big(\f{1323}{100}+\ls{-15}p\big)p
\mod {p^2}&\t{if $p\e 3,5,6\mod 7$,}
\endcases
\\&\Ls {-15}p\sum_{k=0}^{p-1}\f{\b{2k}k\b{3k}k\b{6k}{3k}}
{(-15)^{3k}(2k-1)}\\&\q \e \cases
\f{2542}{225}y^2-\f{2138}{1125}p\mod {p^2}&\t{if $p=x^2+7y^2\e
1,2,4\mod 7$,}
\\-\f 8{15}R_7(p)-\f{112}{125}p\mod {p^2}&\t{if $p\e 3,5,6\mod 7$.}
\endcases\endalign$$
\endpro

\pro{Conjecture 2.25} Let $p>3$ be a prime. Then
$$\align&\sum_{k=0}^{p-1}\f{\b{2k}k^2\b{4k}{2k}}{(-12288)^k(k+1)}
\\&\e\cases 16x^2-14p\mod {p^2} &\t{if $p\e 1\mod {12}$
and so $p=x^2+9y^2$,}
\\-8x^2+14p\mod {p^2} &\t{if $p\e 5\mod {12}$ and so
$2p=x^2+9y^2$,}
\\-24\b{[p/3]}{[p/12]}^2\mod p&\t{if $p\e 7\mod {12}$,}
\\48\b{[p/3]}{[p/12]}^2\mod p&\t{if $p\e 11\mod {12}$,}
\endcases
\\&\sum_{k=0}^{p-1}\f{\b{2k}k^2\b{4k}{2k}}{(-12288)^k(k+1)^2}
\\&\e\cases -128x^2+75p\mod {p^2} &\t{if $p\e 1\mod {12}$
and so $p=x^2+9y^2$,}
\\64x^2-77p\mod {p^2} &\t{if $p\e 5\mod {12}$ and so
$2p=x^2+9y^2$,}
\\-176\b{[p/3]}{[p/12]}^2\mod p&\t{if $p\e 7\mod {12}$,}
\\352\b{[p/3]}{[p/12]}^2\mod p&\t{if $p\e 11\mod {12}$,}
\endcases
\\&\sum_{k=0}^{p-1}\f{\b{2k}k^2\b{4k}{2k}}{(-12288)^k(2k-1)}
\\&\e\cases -\f {13}{4}x^2+\f{93}{64}p\mod {p^2} &\t{if $p\e 1\mod {12}$
and so $p=x^2+9y^2$,}
\\\f{13}8x^2-\f{93}{64}p\mod {p^2} &\t{if $p\e 5\mod {12}$ and so
$2p=x^2+9y^2$,}
\\\f 3{16}\b{[p/3]}{[p/12]}^2\mod p&\t{if $p\e 7\mod {12}$,}
\\-\f 3{8}\b{[p/3]}{[p/12]}^2\mod p&\t{if $p\e 11\mod {12}$.}
\endcases\endalign$$
\endpro

\pro{Conjecture 2.26} Let $p$ be a prime with $p\not=2,5$. Then
$$\align &\sum_{k=0}^{p-1}\f{\b{2k}k^2\b{4k}{2k}}{(-1024)^k}
\\&\e\cases 4x^2-2p-\f{p^2}{4x^2}\mod {p^3}&\t{if $p\e
1,9\mod{20}$ and so $p=x^2+5y^2$,}
\\2p-2x^2+\f{p^2}{2x^2}\mod {p^3}&\t{if $p\e
3,7\mod{20}$ and so $2p=x^2+5y^2$,}
\\\f{2p^2}{\b{(p-1)/2}{[p/20]}\b{(p-1)/2}{[3p/20]}}\mod {p^3}
& \t{if $p\e 11\mod {20}$,}
\\\f{2p^2}{9\b{(p-1)/2}{[p/20]}\b{(p-1)/2}{[3p/20]}}\mod {p^3}
& \t{if $p\e 13\mod {20}$,}
\\\f{6p^2}{7\b{(p-1)/2}{[p/20]}\b{(p-1)/2}{[3p/20]}}
\mod {p^3}& \t{if $p\e 17\mod {20}$,}
\\\f{2p^2}{21\b{(p-1)/2}{[p/20]}\b{(p-1)/2}{[3p/20]}}\mod {p^3}
& \t{if $p\e 19\mod {20}$.}
\endcases\endalign$$
\endpro
\par{\bf Remark 2.3} For any prime $p\not=2,5$, the congruence for
$\sum_{k=0}^{p-1}\f{\b{2k}k^2\b{4k}{2k}}{(-1024)^k}$ modulo $p^2$
was first conjectured by Z.W. Sun in [Su1].

 \pro{Conjecture 2.27}
Let $p$ be a prime with $p\not=2,5$. Then

$$\sum_{k=0}^{p-1}\f{\b{2k}k^2\b{4k}{2k}}{(-1024)^k(k+1)}
\e\cases \f 85\b{\f{p-1}2}{[\f p{20}]}\b{\f{p-1}2}{[\f
{3p}{20}]}\mod p& \t{if $p\e 1,3,7,9\mod {20}$,}
\\\f{4}{5}\b{\f{p-1}2}{[\f p{20}]}\b{\f{p-1}2}{[\f
{3p}{20}]}\mod p& \t{if $p\e 11\mod {20}$,}
\\\f{36}{5}\b{\f{p-1}2}{[\f p{20}]}\b{\f{p-1}2}{[\f
{3p}{20}]}\mod p& \t{if $p\e 13\mod {20}$,}
\\\f{28}{15}\b{\f{p-1}2}{[\f p{20}]}\b{\f{p-1}2}{[\f
{3p}{20}]}\mod p& \t{if $p\e 17\mod {20}$,}
\\\f{84}{5}\b{\f{p-1}2}{[\f p{20}]}\b{\f{p-1}2}{[\f
{3p}{20}]}\mod p& \t{if $p\e 19\mod {20}$.}
\endcases$$
Moreover, if $\sls{-5}p=1$, then
$$\align&\sum_{k=0}^{p-1}\f{\b{2k}k^2\b{4k}{2k}}{(-1024)^k(k+1)}
\\&\e\cases \f{32}{5}x^2-\f{22}{5}p\mod {p^2}&\t{if $p\e
1,9\mod{20}$ and so $p=x^2+5y^2$,}
\\-\f{16}{5}x^2+\f{22}{5}p\mod {p^2}&\t{if $p\e
3,7\mod{20}$ and so $2p=x^2+5y^2$,}\endcases
\\&\sum_{k=0}^{p-1}\f{\b{2k}k^2\b{4k}{2k}}{(-1024)^k(k+1)^2}
\\&\e\cases -\f{32}{5}x^2+\f{7}{5}p\mod {p^2}&\t{if $p\e
1,9\mod{20}$ and so $p=x^2+5y^2$,}
\\\f{16}{5}x^2-\f{17}{5}p\mod {p^2}&\t{if $p\e
3,7\mod{20}$ and so $2p=x^2+5y^2$,}\endcases
\\&\sum_{k=0}^{p-1}\f{\b{2k}k^2\b{4k}{2k}}{(-1024)^k(2k-1)}
\\&\e\cases -\f{31}{10}x^2+\f{103}{80}p\mod {p^2}&\t{if $p\e
1,9\mod{20}$ and so $p=x^2+5y^2$,}
\\\f{31}{20}x^2-\f{103}{80}p\mod {p^2}&\t{if $p\e
3,7\mod{20}$ and so $2p=x^2+5y^2$;}\endcases
\endalign$$
if $\sls{-5}p=-1$, then
$$\align
&\sum_{k=0}^{p-1}\f{\b{2k}k^2\b{4k}{2k}}{(-1024)^k(k+1)^2} \e \f
{22}3\sum_{k=0}^{p-1}\f{\b{2k}k^2\b{4k}{2k}}{(-1024)^k(k+1)}+
\big(8(-1)^{\f{p-1}2}-1\big)p\mod {p^2},
\\&\sum_{k=0}^{p-1}\f{\b{2k}k^2\b{4k}{2k}}{(-1024)^k(2k-1)}
\e -\f
3{32}\sum_{k=0}^{p-1}\f{\b{2k}k^2\b{4k}{2k}}{(-1024)^k(k+1)}-\f
38(-1)^{\f{p-1}2}p\mod {p^2}.\endalign$$
\endpro

\pro{Conjecture 2.28} Let $p$ be a prime with $\sls{-6}p=-1$. Then
$$\sum_{k=0}^{p-1}\f{\b{2k}k^2\b{3k}k}{216^k}\e\cases
-\f 75\cdot \f{p^2}{\b{(p-1)/2}{[p/24]}\b{(p-1)/2}{[5p/24]}} \mod
{p^3}& \t{if $p\e 13\mod {24}$,}
\\\f 75\cdot \f{p^2}{\b{(p-1)/2}{[p/24]}\b{(p-1)/2}{[5p/24]}} \mod
{p^3}& \t{if $p\e 17\mod {24}$,}
\\\f 1{11}\cdot \f{p^2}{\b{(p-1)/2}{[p/24]}\b{(p-1)/2}{[5p/24]}} \mod
{p^3}& \t{if $p\e 19\mod {24}$,}
\\-\f 1{11}\cdot \f{p^2}{\b{(p-1)/2}{[p/24]}\b{(p-1)/2}{[5p/24]}} \mod
{p^3}& \t{if $p\e 23\mod {24}$}
\endcases$$
and
$$\sum_{k=0}^{p-1}\f{\b{2k}k^2\b{4k}{2k}}{48^{2k}}\e\cases
\f 15\cdot \f{p^2}{\b{(p-1)/2}{[p/24]}\b{(p-1)/2}{[5p/24]}} \mod
{p^3}& \t{if $p\e 13,17\mod {24}$,}
\\-\f 1{77}\cdot \f{p^2}{\b{(p-1)/2}{[p/24]}\b{(p-1)/2}{[5p/24]}} \mod
{p^3}& \t{if $p\e 19,23\mod {24}$.}
\endcases$$
\endpro

 \pro{Conjecture 2.29}
Let $p>3$ be a prime. Then
$$\sum_{k=0}^{p-1}\f{\b{2k}k^2\b{3k}k}{216^k(k+1)}\e\cases
\f{7}{8}\b{(p-1)/2}{[p/24]}\b{(p-1)/2}{[5p/24]} \mod p& \t{if $p\e
1,11\mod {24}$,}
\\-\f{7}{8}\b{(p-1)/2}{[p/24]}\b{(p-1)/2}{[5p/24]} \mod p& \t{if
$p\e 5,7\mod {24}$,}
\\ -\f 5{8}\b{(p-1)/2}{[p/24]}\b{(p-1)/2}{[5p/24]}
\mod p& \t{if $p\e 13\mod {24}$,}
\\\f 5{8}\b{(p-1)/2}{[p/24]}\b{(p-1)/2}{[5p/24]} \mod
p& \t{if $p\e 17\mod {24}$,}
\\\f {77}{8}\b{(p-1)/2}{[p/24]}\b{(p-1)/2}{[5p/24]} \mod
p& \t{if $p\e 19\mod {24}$,}
\\-\f {77}{8}\b{(p-1)/2}{[p/24]}\b{(p-1)/2}{[5p/24]} \mod
p& \t{if $p\e 23\mod {24}$.}
\endcases$$
Moreover, if $\sls{-6}p=1$, then
$$\align
&\sum_{k=0}^{p-1}\f{\b{2k}k^2\b{3k}k}{216^k(k+1)}\e\cases \f
72x^2-\f 32p\mod {p^2}&\t{if $p=x^2+6y^2\e 1,7\mod
{24}$,}\\7x^2-2p\mod {p^2}&\t{if $p=2x^2+3y^2\e 5,11\mod
{24}$,}\endcases
\\&\sum_{k=0}^{p-1}\f{\b{2k}k^2\b{3k}k}{216^k(k+1)^2}\e\cases \f
{43}4x^2-\f {25}4p\mod {p^2}&\t{if $p=x^2+6y^2\e 1,7\mod
{24}$,}\\\f{43}2x^2-\f{13}2p\mod {p^2}&\t{if $p=2x^2+3y^2\e 5,11\mod
{24}$,}\endcases
\\&\sum_{k=0}^{p-1}\f{\b{2k}k^2\b{3k}k}{216^k(2k-1)}\e\cases
-\f{23}9x^2+\f 76p\mod {p^2}&\t{if $p=x^2+6y^2\e 1,7\mod
{24}$,}\\-\f{46}9x^2+\f{25}{18}p\mod {p^2}&\t{if $p=2x^2+3y^2\e
5,11\mod {24}$;}\endcases
\endalign$$
if $\sls{-6}p=-1$, then
$$\align
&\sum_{k=0}^{p-1}\f{\b{2k}k^2\b{3k}k}{216^k(k+1)^2}\e \f{13}2
\sum_{k=0}^{p-1}\f{\b{2k}k^2\b{3k}k}{216^k(k+1)}-\Big(1+\f 32\Ls
p3\Big)p\mod {p^2},
\\&\sum_{k=0}^{p-1}\f{\b{2k}k^2\b{3k}k}{216^k(2k-1)}\e\f 29
\sum_{k=0}^{p-1}\f{\b{2k}k^2\b{3k}k}{216^k(k+1)}-\f 16\Ls p3p\mod
{p^2}.\endalign$$

\pro{Conjecture 2.30} Let $p>3$ be a prime. Then
$$\sum_{k=0}^{p-1}\f{\b{2k}k^2\b{4k}{2k}}{48^{2k}(k+1)}\e\cases
\f{1}{3}\b{(p-1)/2}{[p/24]}\b{(p-1)/2}{[5p/24]} \mod p& \t{if $p\e
1,5\mod {24}$,}
\\-\f{1}{3}\b{(p-1)/2}{[p/24]}\b{(p-1)/2}{[5p/24]} \mod p& \t{if
$p\e 7,11\mod {24}$,}
\\ \f 5{3}\b{(p-1)/2}{[p/24]}\b{(p-1)/2}{[5p/24]}
\mod p& \t{if $p\e 13,17\mod {24}$,}
\\-\f {77}3\b{(p-1)/2}{[p/24]}\b{(p-1)/2}{[5p/24]} \mod
p& \t{if $p\e 19,23\mod {24}$.}
\endcases$$
Moreover,
if $\sls{-6}p=1$, then
$$\align&\sum_{k=0}^{p-1}\f{\b{2k}k^2\b{4k}{2k}}{48^{2k}(k+1)}\e\cases
\f 43x^2+\f 23p\mod {p^2}&\t{if $p=x^2+6y^2\e 1,7\mod {24}$,}\\-\f
83x^2-\f 23p\mod {p^2}&\t{if $p=2x^2+3y^2\e 5,11\mod
{24}$},\endcases
\\&\sum_{k=0}^{p-1}\f{\b{2k}k^2\b{4k}{2k}}{48^{2k}(k+1)^2}\e\cases \f
{280}9x^2-\f {157}9p\mod {p^2}&\t{if $p=x^2+6y^2\e 1,7\mod
{24}$,}\\-\f {560}9x^2+\f{139}9p\mod {p^2}&\t{if $p=2x^2+3y^2\e
5,11\mod {24}$},\endcases
\\&\sum_{k=0}^{p-1}\f{\b{2k}k^2\b{4k}{2k}}{48^{2k}(2k-1)}\e\cases
-\f{55}{18}x^2+\f {49}{36}p\mod {p^2}&\t{if $p=x^2+6y^2\e 1,7\mod
{24}$,}\\\f{55}9x^2-\f{49}{36}p\mod {p^2}&\t{if $p=2x^2+3y^2\e
5,11\mod {24}$};\endcases\endalign$$ if $\sls{-6}p=-1$, then
$$\align
&\Ls p3\sum_{k=0}^{p-1}\f{\b{2k}k^2\b{4k}{2k}}{48^{2k}(k+1)} \e-\f
83\sum_{k=0}^{p-1}\f{\b{2k}k^2\b{3k}{k}}{216^k(k+1)} +\f p3\Big(2\Ls
p3+4\Big)\mod {p^2},
\\&\Ls p3\sum_{k=0}^{p-1}\f{\b{2k}k^2\b{4k}{2k}}{48^{2k}(k+1)^2}
\e-\f{176}9\sum_{k=0}^{p-1}\f{\b{2k}k^2\b{3k}{k}}{216^k(k+1)} +\f
p9\Big(35\Ls p3-8\Big)\mod {p^2},
 \\&\Ls
p3\sum_{k=0}^{p-1}\f{\b{2k}k^2\b{4k}{2k}}{48^{2k}(2k-1)}
\e-\f{1}9\sum_{k=0}^{p-1}\f{\b{2k}k^2\b{3k}{k}}{216^k(k+1)} +\f
p{36}\Big(\Ls p3-6\Big)\mod {p^2}.\endalign$$
\endpro

\par{\bf Remark 2.4} Let $p$ be a prime with $p>3$. In
[S10], the author conjectured that if $p\e 1,5,7,11\mod{24}$, then
$$\align\sum_{k=0}^{p-1}\f{\b{2k}k^2\b{3k}k}{216^k}
&\e\Ls p3\sum_{k=0}^{p-1}\f{\b{2k}k^2\b{4k}{2k}}{48^{2k}}
\\&
\e\cases 4x^2-2p-\f{p^2}{4x^2}\mod{p^3}&\t{if $p\e 1,7\mod {24}$ and
so $p=x^2+6y^2$,}
\\8x^2-2p-\f{p^2}{8x^2}\mod {p^3}
&\t{if $p\e 5,11\mod {24}$ and so $p=2x^2+3y^2$;}
\endcases\endalign$$ if $p\e 13,17,19,23\mod{24}$, then
$$\sum_{k=0}^{p-1}\f{\b{2k}k^2\b{3k}k}{216^k}
\e -7\Ls p3\sum_{k=0}^{p-1}\f{\b{2k}k^2\b{4k}{2k}}{48^{2k}} \mod
{p^3}.$$ The corresponding congruences modulo $p^2$ were conjectured
by Z.W. Sun in [Su1].

\pro{Conjecture 2.31} Let $p>5$ be a prime.
\par $(\t{\rm i})$ If
$p\e 1,19\mod{30}$ and so $p=x^2+15y^2$, then
$$\align
&\sum_{k=0}^{p-1}\f{\b{2k}k^2\b{3k}{k}} {(-27)^k(k+1)}\e
-84y^2+2p\mod{p^2},
\\&\sum_{k=0}^{p-1}\f{\b{2k}k^2\b{3k}{k}} {(-27)^k(k+1)^2}\e
-96y^2\mod{p^2},
\\&\sum_{k=0}^{p-1}\f{\b{2k}k^2\b{3k}{k}} {(-27)^k(2k-1)}\e
\f{148}3y^2-\f{26}9p\mod{p^2},
\\& \sum_{k=0}^{p-1}\f{\b{2k}k^2\b{3k}{k}}
{15^{3k}(k+1)}\e 60y^2+2p \mod{p^2},
\\& \sum_{k=0}^{p-1}\f{\b{2k}k^2\b{3k}{k}}
{15^{3k}(k+1)^2}\e -1320y^2+36p \mod{p^2},
\\&
\sum_{k=0}^{p-1}\f{\b{2k}k^2\b{3k}{k}} {15^{3k}(2k-1)}\e
\f{3508}{75}y^2-\f{1954}{1125}p
 \mod{p^2}.\endalign$$
\par $(\t{\rm ii})$ If $p\e
17,23\mod{30}$ and so $p=3x^2+5y^2$, then
$$\align
&\sum_{k=0}^{p-1}\f{\b{2k}k^2\b{3k}{k}} {(-27)^k(k+1)}\e
28y^2-2p\mod{p^2},
\\&\sum_{k=0}^{p-1}\f{\b{2k}k^2\b{3k}{k}} {(-27)^k(k+1)^2}\e
32y^2-2p\mod{p^2},
\\&\sum_{k=0}^{p-1}\f{\b{2k}k^2\b{3k}{k}} {(-27)^k(2k-1)}\e
-\f{148}9y^2+\f{26}9p\mod{p^2},
\\& \sum_{k=0}^{p-1}\f{\b{2k}k^2\b{3k}{k}}
{15^{3k}(k+1)}\e -20y^2-2p \mod{p^2},
\\& \sum_{k=0}^{p-1}\f{\b{2k}k^2\b{3k}{k}}
{15^{3k}(k+1)^2}\e 440y^2-38p \mod{p^2},
\\&
\sum_{k=0}^{p-1}\f{\b{2k}k^2\b{3k}{k}} {15^{3k}(2k-1)}\e
-\f{3508}{225}y^2+\f{1954}{1125}p
 \mod{p^2}.\endalign$$
\par $(\t{\rm iii})$ If $\sls{-15}p=-1$, then
$$\sum_{k=0}^{p-1}\f{\b{2k}k^2\b{3k}{k}}
{(-27)^k(k+1)}\e\cases \f 25\cdot 5^{-[\f
p3]}\b{[p/3]}{[p/15]}^2\mod p&\t{if $p\e 7\mod {30}$,}
\\\f  1{10}\cdot 5^{-[\f p3]}\b{[p/3]}{[p/15]}^2\mod p&\t{if $p\e 11\mod
{30}$,}
\\\f{32}{5}\cdot 5^{-[\f p3]}\b{[p/3]}{[p/15]}^2\mod p&\t{if $p\e 13\mod
{30}$,}
\\\f{8}{5}\cdot 5^{-[\f p3]}\b{[p/3]}{[p/15]}^2\mod p&\t{if $p\e 29\mod
{30}$.}
\endcases$$
Moreover,
$$\align &\sum_{k=0}^{p-1}\f{\b{2k}k^2\b{3k}{k}}
{(-27)^k(k+1)^2}\e\f{13}2 \sum_{k=0}^{p-1}\f{\b{2k}k^2\b{3k}{k}}
{(-27)^k(k+1)}+\Big(3\Ls p3-1\Big)p\mod {p^2},
\\&\sum_{k=0}^{p-1}\f{\b{2k}k^2\b{3k}{k}}
{(-27)^k(2k-1)}\e-\f{16}9 \sum_{k=0}^{p-1}\f{\b{2k}k^2\b{3k}{k}}
{(-27)^k(k+1)}-\f 83\Ls p3p\mod {p^2},
\\&\sum_{k=0}^{p-1}\f{\b{2k}k^2\b{3k}{k}}
{15^{3k}(k+1)}\e-20 \sum_{k=0}^{p-1}\f{\b{2k}k^2\b{3k}{k}}
{(-27)^k(k+1)}-12\Ls p3p\mod {p^2},
\\&\sum_{k=0}^{p-1}\f{\b{2k}k^2\b{3k}{k}}
{15^{3k}(k+1)^2}\e-130 \sum_{k=0}^{p-1}\f{\b{2k}k^2\b{3k}{k}}
{(-27)^k(k+1)}-\Big(1+111\Ls p3\Big)p\mod {p^2},
\\&\sum_{k=0}^{p-1}\f{\b{2k}k^2\b{3k}{k}}
{15^{3k}(2k-1)}\e-\f{64}{225} \sum_{k=0}^{p-1}\f{\b{2k}k^2\b{3k}{k}}
{(-27)^k(k+1)}-\f{152}{375}\Ls p3p\mod {p^2}.
\endalign$$
\endpro

\pro{Conjecture 2.32} Let $p$ be a prime with $p>5$ and
$$t_p=\f{\b{(p-1)/2}{[7p/40]}\b{(p-1)/2}{[9p/40]}
\b{[3p/40]}{[p/40]}}{\b{[19p/40]}{[p/20]}}.$$
 \par $(\t{\rm i})$ If
$\sls{-10}p=1$, then
$$\sum_{k=0}^{p-1}\f{\b{2k}k^2\b{4k}{2k}}{12^{4k}(k+1)}
\e\cases -\f{49}{15}t_p\mod p&\t{if $p\e 1,23\mod {40}$,}
\\-\f{49}{5}t_p\mod p&\t{if $p\e 7\mod {40}$,}
\\\f{7}{5}t_p\mod p&\t{if $p\e 9\mod {40}$,}
\\\f{833}{195}t_p\mod p&\t{if $p\e 11\mod {40}$,}
\\-\f{833}{285}t_p\mod p&\t{if $p\e 13\mod {40}$,}
\\-\f{98}{555}t_p\mod p&\t{if $p\e 19\mod {40}$,}
\\-\f{98}{55}t_p\mod p&\t{if $p\e 37\mod {40}$.}
\endcases$$
Moreover,
$$\align&\sum_{k=0}^{p-1}\f{\b{2k}k^2\b{4k}{2k}}{12^{4k}(k+1)}
\e\cases-\f{196}{15}x^2+\f{226}{15}p\mod {p^2}&\t{if $p=x^2+10y^2\e
1,9,11,19\mod {40}$,}
\\\f{392}{15}x^2-\f{226}{15}p\mod {p^2}&\t{if $p=2x^2+5y^2\e 7,13,19,37\mod {40}$,}
\endcases
\\&\sum_{k=0}^{p-1}\f{\b{2k}k^2\b{4k}{2k}}{12^{4k}(k+1)^2}
\\&\q\e\cases\f{10088}{45}x^2-\f{6113}{45}p\mod {p^2}&\t{if
$p=x^2+10y^2\e 1,9,11,19\mod {40}$,}
\\-\f{20176}{45}x^2+\f{6023}{45}p\mod {p^2}&\t{if $p=2x^2+5y^2\e 7,13,19,37\mod {40}$,}
\endcases
\\&\sum_{k=0}^{p-1}\f{\b{2k}k^2\b{4k}{2k}}{12^{4k}(2k-1)}
\\&\q\e\cases -\f{883}{270}x^2+\f{2393}{1620}p\mod {p^2}&\t{if
$p=x^2+10y^2\e 1,9,11,19\mod {40}$,}
\\\f{883}{135}x^2-\f{2393}{1620}p\mod {p^2}&\t{if $p=2x^2+5y^2\e 7,13,23,37\mod
{40}$}\endcases\endalign$$ and
$$\align &\sum_{k=0}^{p-1}\f{\b{2k}k^2\b{4k}{2k}}{12^{4k}}
\\&\e\cases 4x^2-2p-\f{p^2}{4x^2}\mod {p^3}&\t{if $p\e 1,9,11,19\mod
{40}$ and so $p=x^2+10y^2$,}
\\2p-8x^2+\f{p^2}{8x^2}\mod {p^3}&\t{if $p\e 7,13,23,37\mod
{40}$ and so $p=2x^2+5y^2$.}\endcases\endalign$$
\par $(\t{\rm ii})$ If
$\sls{-10}p=-1$, then
$$\sum_{k=0}^{p-1}\f{\b{2k}k^2\b{4k}{2k}}{12^{4k}(k+1)}
\e\cases -\f{21}{5}t_p\mod p&\t{if $p\e 3\mod {40}$,}
\\-\f{4446}{155}t_p\mod p&\t{if $p\e 17\mod {40}$,}
\\-\f{189}{5}t_p\mod p&\t{if $p\e 21\mod {40}$,}
\\-\f{702}{5}t_p\mod p&\t{if $p\e 27\mod {40}$,}
\\\f{66}{5}t_p\mod p&\t{if $p\e 29\mod {40}$,}
\\\f{1026}{5}t_p\mod p&\t{if $p\e 31\mod {40}$,}
\\-\f{462}{5}t_p\mod p&\t{if $p\e 33\mod {40}$,}
\\-\f{858}{85}t_p\mod p&\t{if $p\e 39\mod {40}.$}
\endcases$$
Moreover,
$$\align&\sum_{k=0}^{p-1}\f{\b{2k}k^2\b{4k}{2k}}{12^{4k}(k+1)^2}
\e \f{22}3\sum_{k=0}^{p-1}\f{\b{2k}k^2\b{4k}{2k}}{12^{4k}(k+1)}
+\Big(\f{256}3\Ls p5-1\Big)p\mod {p^2},
\\&\sum_{k=0}^{p-1}\f{\b{2k}k^2\b{4k}{2k}}{12^{4k}(2k-1)}
\e \f 1{216}\sum_{k=0}^{p-1}\f{\b{2k}k^2\b{4k}{2k}}{12^{4k}(k+1)}
+\f{16}{81}\Ls p5p\mod {p^2}\endalign$$ and
$$\Big(\sum_{k=0}^{p-1}\f{\b{2k}k^2\b{4k}{2k}}{12^{4k}(k+1)}\Big)
\Big(\sum_{k=0}^{p-1}\f{\b{2k}k^2\b{4k}{2k}}{12^{4k}}\Big) \e
-\f{49}{15}p^2\mod {p^3}.$$
\endpro
\par{\bf Remark 2.5} Let $p>3$ be a prime. In [Su1], Z.W. Sun
conjectured the congruence for
$\sum_{k=0}^{p-1}\f{\b{2k}k^2\b{4k}{2k}}{12^{4k}}\mod {p^2}$.
 \pro{Conjecture 2.33} Let $p$ be a
prime with $p\not=2,11$ and
$r_p=\b{[3p/11]}{[p/11]}\b{[6p/11]}{[3p/11]}/\b{[4p/11]}{[2p/11]}$.
\par $(\t{\rm i})$ If $p\e 1,3,4,5,9\mod {11}$
and so $4p=x^2+11y^2$, then
$$\sum_{k=0}^{p-1}\f{\b{2k}k^2\b{3k}{k}}
{64^k(k+1)}\e\cases \f{25}{22}r_p^2\mod p&\t{if $p\e
1,4,5,9\mod{11}$,}
\\\f 2{11}r_p^2\mod p&\t{if $p\e 3\mod {11}$.}
\endcases$$
Moreover,
$$\align&\sum_{k=0}^{p-1}\f{\b{2k}k^2\b{3k}{k}}
{64^k(k+1)}\e -\f{25}2y^2+2p\mod {p^2},
\\&\sum_{k=0}^{p-1}\f{\b{2k}k^2\b{3k}{k}}
{64^k(k+1)^2}\e -\f{83}4y^2+2p\mod {p^2},
\\&\sum_{k=0}^{p-1}\f{\b{2k}k^2\b{3k}{k}}
{64^k(2k-1)}\e  \f {23}4y^2-\f 78p\mod {p^2},
\\&\Ls{-2}p\sum_{k=0}^{p-1}\f{\b{2k}k\b{3k}k\b{6k}{3k}}{(-32)^{3k}(k+1)}
\e -26y^2+2p\mod {p^2},
\\&\Ls{-2}p\sum_{k=0}^{p-1}\f{\b{2k}k\b{3k}k\b{6k}{3k}}{(-32)^{3k}(k+1)^2}
\e 148y^2-\Big(24+\Ls{-2}p\Big)p\mod {p^2},
\\&\Ls{-2}p\sum_{k=0}^{p-1}\f{\b{2k}k\b{3k}k\b{6k}{3k}}{(-32)^{3k}(2k-1)}
\e \f{73}{8}y^2-\f{467}{256}p\mod {p^2}.
\endalign$$
\par $(\t{\rm ii})$ If $p\e 2,6,7,8,10\mod {11}$, then
$$ \sum_{k=0}^{p-1}\f{\b{2k}k^2\b{3k}{k}}
{64^k(k+1)}\e \cases -\f{50}{11}r_p^2\mod p&\t{if $p\e 2\mod {11}$,}
\\-\f{32}{11}r_p^2\mod p&\t{if $p\e 6\mod{11}$,}
\\-\f{2}{11}r_p^2\mod p&\t{if $p\e 7\mod{11}$,}
\\-\f{72}{11}r_p^2\mod p&\t{if $p\e 8\mod{11}$,}
\\-\f{18}{11}r_p^2\mod p&\t{if $p\e 10\mod{11}$.}
\endcases$$
Moreover,
$$\align &\sum_{k=0}^{p-1}\f{\b{2k}k^2\b{3k}{k}}
{64^k(k+1)^2}\e\f{13}2 \sum_{k=0}^{p-1}\f{\b{2k}k^2\b{3k}{k}}
{64^k(k+1)}\mod {p^2},
\\&\sum_{k=0}^{p-1}\f{\b{2k}k^2\b{3k}{k}} {64^k(2k-1)}
\e \f 34\sum_{k=0}^{p-1}\f{\b{2k}k^2\b{3k}{k}} {64^k(k+1)}+\f
38p\mod {p^2},
\\&\Ls{-2}p\sum_{k=0}^{p-1}\f{\b{2k}k\b{3k}k\b{6k}{3k}}
{(-32)^{3k}(k+1)}\e
\f{128}{15}\sum_{k=0}^{p-1}\f{\b{2k}k^2\b{3k}{k}} {64^k(k+1)} -\f
25p\mod {p^2},
\\&\Ls{-2}p\sum_{k=0}^{p-1}\f{\b{2k}k\b{3k}k\b{6k}{3k}}
{(-32)^{3k}(k+1)^2}\e
\f{5888}{75}\sum_{k=0}^{p-1}\f{\b{2k}k^2\b{3k}{k}} {64^k(k+1)}
+\Big(\f{608}{25}-\Ls{-2}p\Big)p\mod {p^2},
\\&\Ls{-2}p\sum_{k=0}^{p-1}\f{\b{2k}k\b{3k}k\b{6k}{3k}}
{(-32)^{3k}(2k-1)}\e -\f 18\sum_{k=0}^{p-1}\f{\b{2k}k^2\b{3k}{k}}
{64^k(k+1)}-\f{51}{256}p\mod {p^2}\endalign$$ and
$$\Big(\sum_{k=0}^{p-1}\f{\b{2k}k^2\b{3k}{k}}
{64^k(k+1)}\Big)\Big(\sum_{k=0}^{p-1}\f{\b{2k}k^2\b{3k}{k}}
{64^k}\Big)\e \f{25}{22}p^2\mod{p^3}.$$
\endpro
\par{\bf Remark 2.6} Let $p$ be a prime with $p\not=2,3,11,13$. In
[S10] the author conjectured that if $p\e 1,3,4,5,9\mod{11}$ and so
$4p=x^2+11y^2$, then
$$\sum_{k=0}^{p-1}\f{\b{2k}k^2\b{3k}k}{64^k}\e
\Ls{-2}p\sum_{k=0}^{p-1}\f{\b{2k}k\b{3k}k\b{6k}{3k}}{(-32)^{3k}} \e
x^2-2p-\f {p^2}{x^2}\mod {p^3};$$ if $p\e 2,6,7,8,10\mod{11}$ and
$f=[\f p{11}]$, then
$$\align \sum_{k=0}^{p-1}\f{\b{2k}k^2\b{3k}k}{64^k}
&\e
\f{160}{39}\Ls{-2}p\sum_{k=0}^{p-1}\f{\b{2k}k\b{3k}k\b{6k}{3k}}{(-32)^{3k}}
\\&\e\cases -p^2\Ls{5\b{4f}{2f}}{2\b{3f}f\b{6f}{3f}}^2\mod {p^3}&\t{if
$p\e 2\mod {11}$,}
\\-p^2\Ls{13\b{4f}{2f}}{30\b{3f}f\b{6f}{3f}}^2\mod
{p^3}&\t{if $p\e 6\mod {11}$,}
\\-p^2\Ls{85\b{4f}{2f}}{558\b{3f}f\b{6f}{3f}}^2\mod{p^3}&\t{if $p\e
7\mod {11}$,}
\\-p^2\Ls{7\b{4f}{2f}}{148\b{3f}f\b{6f}{3f}}^2\mod{p^3}&\t{if $p\e
8\mod {11}$,}
\\-p^2\Ls{29\b{4f}{2f}}{756\b{3f}f\b{6f}{3f}}^2\mod{p^3}&\t{if $p\e
10\mod {11}$.}
\endcases\endalign$$
The corresponding congruence modulo $p^2$ was conjectured by Z.W.
Sun [Su3].
 \pro{Conjecture 2.34} Let $p>3$ be a prime and
$s_p=\b{[8p/19]}{[p/19]}\b{[10p/19]}{[4p/19]}/\b{[5p/19]}{[2p/19]}$.
\par $(\t{\rm i})$ If $\sls {p}{19}=1$ and so $4p=x^2+19y^2$, then
$$\Ls{-6}p\sum_{k=0}^{p-1}\f{\b{2k}k\b{3k}{k}\b{6k}{3k}}
{(-96)^{3k}(2k-1)} \e\cases -\f{1183}{1368}s_p^2\mod p&\t{if $p\e
1,7,11\mod {19}$,} \\ -\f{1183}{342}s_p^2\mod p&\t{if $p\e 4\mod
{19}$,}
\\ -\f{1183}{2432}s_p^2\mod p&\t{if $p\e 5,6\mod
{19}$,}
\\ -\f{1183}{18392}s_p^2\mod p&\t{if $p\e 9\mod
{19}$,}
\\ -\f{1183}{27702}s_p^2\mod p&\t{if $p\e 16\mod
{19}$,}
\\ -\f{57967}{12312}s_p^2\mod p&\t{if $p\e 17\mod
{19}$.}
\endcases$$
Moreover,
$$\align &\Ls{-6}p\sum_{k=0}^{p-1}\f{\b{2k}k\b{3k}{k}\b{6k}{3k}}
{(-96)^{3k}(2k-1)} \e \f{1183}{72}y^2-\f{4273}{2304}p\mod{p^2},
\\&\Ls{-6}p\sum_{k=0}^{p-1}\f{\b{2k}k\b{3k}{k}\b{6k}{3k}}
{(-96)^{3k}(k+1)} \e -394y^2+2p\mod{p^2},
\\&\Ls{-6}p\sum_{k=0}^{p-1}\f{\b{2k}k\b{3k}{k}\b{6k}{3k}}
{(-96)^{3k}(k+1)^2} \e 6772y^2-\Big(536+\Ls{-6}p\Big)p\mod{p^2}
.\endalign$$
\par $(\t{\rm ii})$ If $\sls {p}{19}=-1$,
 then
$$ \Ls{-6}p\sum_{k=0}^{p-1}\f{\b{2k}k\b{3k}{k}\b{6k}{3k}}
{(-96)^{3k}(2k-1)} \e\cases \f{49}{228}s_p^2\mod p&\t{if $p\e 2\mod
{19}$,}
\\ \f{1}{228}s_p^2\mod p&\t{if $p\e
3\mod {19}$,}
\\ \f{27}{14896}s_p^2\mod p&\t{if $p\e
8\mod {19}$,}
\\ \f{3}{19}s_p^2\mod p&\t{if $p\e
10\mod {19}$,}
\\ \f{121}{57}s_p^2\mod p&\t{if $p\e
12\mod {19}$,}
\\ \f{4}{57}s_p^2\mod p&\t{if $p\e
13\mod {19}$,}
\\ \f{121}{228}s_p^2\mod p&\t{if $p\e
14\mod {19}$,}
\\ \f{27}{76}s_p^2\mod p&\t{if $p\e
15\mod {19}$,}
\\ \f{49}{57}s_p^2\mod p&\t{if $p\e
18\mod {19}$.}
\endcases$$
Moreover,
$$\align &\sum_{k=0}^{p-1}\f{\b{2k}k\b{3k}{k}\b{6k}{3k}}
{(-96)^{3k}(k+1)}\e
-\f{9216}5\sum_{k=0}^{p-1}\f{\b{2k}k\b{3k}{k}\b{6k}{3k}}
{(-96)^{3k}(2k-1)}-270\Ls{-6}pp\mod {p^2},
\\&\sum_{k=0}^{p-1}\f{\b{2k}k\b{3k}{k}\b{6k}{3k}}
{(-96)^{3k}(k+1)^2}\e
\f{46}5\sum_{k=0}^{p-1}\f{\b{2k}k\b{3k}{k}\b{6k}{3k}}
{(-96)^{3k}(k+1)}+\Big(540\Ls{-6}p-1\Big)p\mod {p^2}
\endalign$$
and
$$\Big(\sum_{k=0}^{p-1}\f{\b{2k}k\b{3k}{k}\b{6k}{3k}}
{(-96)^{3k}(2k-1)}\Big)\Big(
\sum_{k=0}^{p-1}\f{\b{2k}k\b{3k}{k}\b{6k}{3k}} {(-96)^{3k}}\Big)\e
-\f{985}{87552}p^2\mod {p^3}.$$
\endpro

 \pro{Conjecture 2.35} Let $p>5$ be a prime. If $\sls
{p}{43}=1$ and so $4p=x^2+43y^2$, then
$$\Ls{-15}p\sum_{k=0}^{p-1}\f{\b{2k}k\b{3k}{k}\b{6k}{3k}}
{(-960)^{3k}(2k-1)} \e
\f{140501}{3600}y^2-\f{4384321}{2304000}p\mod{p^2}.$$ If $\sls
{p}{43}=-1$, then
$$\Big(\sum_{k=0}^{p-1}\f{\b{2k}k\b{3k}{k}\b{6k}{3k}}
{(-960)^{3k}(2k-1)}\Big)\Big(
\sum_{k=0}^{p-1}\f{\b{2k}k\b{3k}{k}\b{6k}{3k}} {(-960)^{3k}}\Big)\e
-\f{933889}{198144000}p^2\mod {p^3}.$$
\endpro

\pro{Conjecture 2.36} Let $p>5$ be a prime.  If $\sls {p}{67}=-1$,
then
$$\Big(\sum_{k=0}^{p-1}\f{\b{2k}k\b{3k}{k}\b{6k}{3k}}
{(-5280)^{3k}(2k-1)}\Big)\Big(
\sum_{k=0}^{p-1}\f{\b{2k}k\b{3k}{k}\b{6k}{3k}} {(-5280)^{3k}}\Big)\e
-\f{155357161}{51365952000}p^2\mod {p^3}.$$
\endpro
\par{\bf Remark 2.7} Suppose that $p>3$ is a prime. In [Su3],
 Z.W. Sun conjectured that
$$\align&\sum_{k=0}^{p-1}\f{\b{2k}k\b{3k}k\b{6k}{3k}}{(-32)^{3k}}
\e\cases 0\mod {p^2}&\t{if $\sls p{11}=-1$,}
\\\ls{-2}p(x^2-2p)\mod {p^2}&\t{if $\sls p{11}=1$ and so $4p=x^2+11y^2$.}
\endcases
\\&\sum_{k=0}^{p-1}\f{\b{2k}k\b{3k}k\b{6k}{3k}}{(-96)^{3k}}
\e\cases 0\mod {p^2}&\t{if $\sls p{19}=-1$,}
\\\ls{-6}p(x^2-2p)\mod {p^2}&\t{if $\sls p{19}=1$ and so $4p=x^2+19y^2$,}
\endcases
\\&
\sum_{k=0}^{p-1}\f{\b{2k}k\b{3k}k\b{6k}{3k}}{(-960)^{3k}} \e\cases
0\mod {p^2}&\t{if $\sls p{43}=-1$ and $p\not=5$,}
\\\ls p{15}(x^2-2p)\mod {p^2}&\t{if $\sls p{43}=1$ and so $4p=x^2+43y^2$,}
\endcases
\\&\sum_{k=0}^{p-1}\f{\b{2k}k\b{3k}k\b{6k}{3k}}{(-5280)^{3k}}
 \e\cases 0\mod {p^2}&\t{if $\sls p{67}=-1$ and $p\not=5,11$,}
\\\ls {-330}p(x^2-2p)\mod {p^2}&\t{if $\sls p{67}=1$ and so $4p=x^2+67y^2$}
\endcases
\endalign$$
and
$$\sum_{k=0}^{p-1}\f{\b{2k}k\b{3k}k\b{6k}{3k}}{(-640320)^{3k}}
\e\cases 0\mod {p^2}&\t{if $\sls p{163}=-1$ and $p\not=5,23,29$,}
\\\ls {-10005}p(x^2-2p)\mod {p^2}&\t{if $\sls p{163}=1$ and so $4p=x^2+163y^2$.}
\endcases$$
These congruences modulo $p$ were proved by the author in [S4]

\pro{Conjecture 2.37} Let $p$ be an odd prime. If $p\e 1\mod 4$ and
so $p=x^2+4y^2$ with $4\mid x-1$, then
$$\sum_{k=0}^{p-1}\f{\b{2k}k^2}{32^k(2k-1)^2}
\e x-\f p{4x}\mod
{p^2}\qtq{and}\sum_{k=0}^{p-1}\f{\b{2k}k^2}{32^k(k+1)^2}\e 8x-7\mod
p.$$ If $p\e 3\mod 4$, then
$$\sum_{k=0}^{p-1}\f{\b{2k}k^2}{32^k(2k-1)^2}\e
\f 12(2p+3-2^{p-1})\b{\f{p-1}2}{\f{p-3}4}\mod {p^2}.$$
\endpro

\pro{Conjecture 2.38} Let $p$ be an odd prime. Then
$$\sum_{k=0}^{p-1}\f{\b{2k}k^2}{(-16)^k(2k-1)^2}
\e \cases (-1)^{\f{p-1}4}\f px\mod{p^2}\q\t{if $p=x^2+4y^2\e 1\mod
4$ and $4\mid x-1$,}
\\2(-1)^{\f{p+1}4}\b{(p-1)/2}{(p-3)/4}\mod p\q\t{if $p\e 3\mod 4$}
\endcases$$
and $$\sum_{k=0}^{p-1}\f{\b{2k}k^2}{(-16)^k(k+1)^2}\e 5\mod
p\qtq{for}p\e 1\mod 4.$$
\endpro

\pro{Conjecture 2.39} Let $p$ be an odd prime. Then
$$\sum_{k=0}^{p-1}\f{\b{2k}k^2}{8^k(2k-1)^2}
\e\cases 2(-1)^{\f{p-1}4}x\mod {p}&\t{if $p=x^2+4y^2\e 1\mod 4$ and
$4\mid x-1$,}
\\2(-1)^{\f{p+1}4}\b{\f{p-1}2}{\f{p-3}4}\mod p&\t{if $p\e 3\mod 4$.}
\endcases$$
\endpro

\par{\bf Remark 2.8} Let $p$ be an odd prime. In [Su4], Z.W. Sun
established the congruences for
$\sum_{k=0}^{p-1}\f{\b{2k}k^2}{m^k(2k-1)}\mod {p^2}$ in the cases
$m=-16,8,32$.

\pro{Conjecture 2.40} Let $a,x\in\Bbb Q$ and $d\in\Bbb Z$ with
$adx\not=0$, where $\Bbb Q$ is the set of rational numbers. Suppose
that $\sum_{k=0}^{p-1}\b ak\b{-1-a}kx^k\e 0\mod p$ for all odd
primes $p$ satisfying $a,x\in\Bbb Z_p$ and $\ls{d}p=-1$. Then there
is a constant $c\in\Bbb Q$ such that for all odd primes $p$ with
$c\in\Bbb Z_p$ and $\ls{d}p=-1$,
$$
\sum_{k=0}^{p-1}\b{2k}k\b ak\b{-1-a}k(x(1-x))^k\e
c\Big(\sum_{k=0}^{p-1}\b ak\b{-1-a}kx^k\Big)^2\mod {p^3}.$$
\endpro
\par {\bf Remark 2.9} By [S6], for any odd prime $p$ and $a,x\in\Bbb Z_p$,
$$\Big(\sum_{k=0}^{p-1}\b ak\b{-1-a}kx^k\Big)^2\e
\sum_{k=0}^{p-1}\b {2k}k\b ak\b{-1-a}k(x(1-x))^k\mod{p^2}.$$
\pro{Conjecture 2.41} Let $a,m\in\Bbb Q$ and $d\in\Bbb Z$ with
$adm\not=0$.  Suppose that $$\sum_{k=0}^{p-1}\f{\b{2k}k\b
ak\b{-1-a}k}{m^k}\e 0\mod {p^2}$$ for all odd primes $p$ satisfying
$a,m\in\Bbb Z_p$, $p\nmid m$ and $\ls{d}p=-1$. Then there is a
constant $c\in\Bbb Q$ such that for all odd primes $p$ with
$c,m\in\Bbb Z_p$, $p\nmid m$ and $\ls{d}p=-1$,
$$\Big(\sum_{k=0}^{p-1}\f {\b{2k}k\b
ak\b{-1-a}k}{(2k-1)m^k}\Big)\Big(\sum_{k=0}^{p-1}\f{\b{2k}k\b
ak\b{-1-a}k}{m^k}\Big)\e cp^2\mod {p^3}.$$ Moreover, there is a
constant $C\in\Bbb Q$ such that for all odd primes $p$ with
$C,m\in\Bbb Z_p$, $p\nmid m$ and $\ls{d}p=1$,
$$\sum_{k=0}^{p-1}\f {\b{2k}k\b
ak\b{-1-a}k}{(2k-1)m^k}\e C\sum_{k=0}^{p-1}\f{\b{2k}k\b
ak\b{-1-a}k}{m^k}\mod p.$$
\endpro

\section*{3. Conjectures for congruences involving Ap\'ery-like numbers}
\par\q
With the help of Maple, we discover the following conjectures
involving Ap\'ery-like numbers.
 \pro{Conjecture 3.1} Let $p$ be an odd prime. Then
$$
\align &\sum_{k=0}^{p-1}\f{\b{2k}kS_k}{16^k(k+1)} \e\cases
4x^2-2p\mod {p^2}&\t{if $p=x^2+4y^2\e 1\mod 4$,}
\\ 0\mod {p^2}&\t{if $p\e 3\mod 4$,}\endcases
\\&\sum_{k=0}^{p-1}\f{\b{2k}kS_k}{16^k(2k-1)} \e\cases 0\mod
{p^2}&\t{if $p\e 1\mod 4$,}
\\-R_1(p)\mod {p^2}&\t{if $p\e 3\mod 4$.}
\endcases\endalign$$\endpro
\pro{Conjecture 3.2} Let $p>3$ be a prime. Then
$$\align &\sum_{k=0}^{p-1}\f{\b{2k}kS_k}{32^k(2k-1)}
\e\cases (-1)^{\f{p-1}2}( \f p2-x^2)\mod {p^2}\qq\t{if $p=x^2+2y^2\e
1,3\mod 8$,}
\\\f {5(-1)^{\f{p-1}2}-4}4R_2(p)\mod {p^2}\q\t{if $p\e 5,7\mod 8$,}
\endcases
\\&\sum_{k=0}^{p-1}\f{\b{2k}kS_k}{800^k(2k-1)}
\\&\q\e\cases -\f{73}{100}(-1)^{\f{p-1}2}( 4x^2-2p)-\f{24}{125}\Ls
3pp\mod {p^2}\ \t{if $p=x^2+2y^2\e 1,3\mod 8$,}
\\\f {45(-1)^{\f{p-1}2}-36}{100}R_2(p)+\f{24}{125}\Ls 3pp
\mod {p^2}\q\t{if $p\e 5,7\mod 8$ and $p\not=5$,}
\endcases
\\&\sum_{k=0}^{p-1}\f{\b{2k}kS_k}{(-768)^k(2k-1)}
\\&\q\e\cases -\f{73}{96}\Ls 3p( 4x^2-2p)-\f{11}{48}(-1)^{\f{p-1}2}p\mod {p^2}\ \t{if $p=x^2+2y^2\e 1,3\mod 8$,}
\\-\f {15-12(-1)^{\f{p-1}2}}{32}\Ls 3pR_2(p)-(-1)^{\f{p-1}2}
\f{11}{48}p \mod {p^2}\q\t{if $p\e 5,7\mod 8$.}
\endcases
\endalign$$
\endpro

\pro{Conjecture 3.3} Let $p$ be an odd prime. Then
$$\align&\sum_{k=0}^{p-1}\f{\b{2k}kS_k}{7^k(2k-1)}
\e\cases \f{124}{49}x^2-\f{46}{49}p\mod {p^2}\qq\t{if $p=x^2+7y^2\e
1,2,4\mod 7$,}
\\\f{64}7\sum_{k=0}^{p-1}\f{\b{2k}k^3}{k+1}+\f{496}{49}p\mod {p^2}
\q\t{if $p\e 3,5,6\mod 7$,}
\endcases
\\&\sum_{k=0}^{p-1}\f{\b{2k}kS_k}{25^k(2k-1)}
\e\cases -\f{124}{175}x^2+\f{326}{875}p\mod {p^2}\qq\t{if
$p=x^2+7y^2\e 1,2,4\mod 7$,}
\\\f{448}{175}\sum_{k=0}^{p-1}\f{\b{2k}k^3}{k+1}+\f{2576}{875}p\mod {p^2}
\q\t{if $p\e 3,5,6\mod 7$ and $p>7$.}
\endcases\endalign$$

\endpro

\pro{Conjecture 3.4} Let $p>3$ be a prime. Then
$$\align&(-1)^{\f{p-1}2}\sum_{k=0}^{p-1}\f{\b{2k}kS_k}{(-32)^k(2k-1)}
\\&\q\e\cases -\f{11}3x^2+\f 76p\mod {p^2}&\t{if $p=x^2+6y^2\e 1,7\mod
{24}$,}
\\\f{22}3x^2-\f 76p\mod {p^2}&\t{if $p=2x^2+3y^2\e 5,11\mod
{24}$,}
\\\f 23\sum_{k=0}^{p-1}\f{\b{2k}k^2\b{3k}k}{216^k(k+1)}+\f p2
\mod {p^2}&\t{if $p\e 13,19\mod {24}$,}
\\-\f 23\sum_{k=0}^{p-1}\f{\b{2k}k^2\b{3k}k}{216^k(k+1)}-\f 56p
\mod {p^2}&\t{if $p\e 17,23\mod {24}$}
\endcases\endalign$$
and
$$\align&(-1)^{\f{p-1}2}\sum_{k=0}^{p-1}\f{\b{2k}kS_k}{64^k(2k-1)}
\\&\q\e\cases -\f{11}6x^2+\f 56p\mod {p^2}&\t{if $p=x^2+6y^2\e 1,7\mod
{24}$,}
\\-\f{11}3x^2+\f 56p\mod {p^2}&\t{if $p=2x^2+3y^2\e 5,11\mod
{24}$,}
\\-\f 13\sum_{k=0}^{p-1}\f{\b{2k}k^2\b{3k}k}{216^k(k+1)}
\mod {p^2}&\t{if $p\e 13,19\mod {24}$,}
\\-\f 13\sum_{k=0}^{p-1}\f{\b{2k}k^2\b{3k}k}{216^k(k+1)}-\f p6
\mod {p^2}&\t{if $p\e 17,23\mod {24}$.}
\endcases\endalign$$

\endpro

\pro{Conjecture 3.5} Let $p$ be a prime with $p>3$. Then
$$\sum_{k=0}^{p-1}\f{\b{2k}kW_k}{(-12)^k(2k-1)}
\e \cases 0\mod {p^2}&\t{if $p\e 1\mod 3$,}
\\-2(2p+1)\b{[2p/3]}{[p/3]}^2\mod {p^2}&\t{if $p\e 2\mod {3}$}.
\endcases$$
\endpro
\pro{Conjecture 3.6} Let $p$ be a prime with $p>3$. Then
$$\Ls p3\sum_{k=0}^{p-1}\f{\b{2k}kW_k}{54^k(2k-1)}
\e \cases -\f{28}9x^2+\f{10}9p\mod {p^2}&\t{if $p=x^2+4y^2\e 1\mod
4$,}\\\f 23R_1(p)-\f 49p\mod {p^2}& \t{if $p\e 3\mod 4$.}
\endcases$$
\endpro

\pro{Conjecture 3.7} Let $p$ be an odd prime. Then
$$\sum_{k=0}^{p-1}\f{\b{2k}kW_k}{8^k(2k-1)} \e \cases -\f{11}2x^2+\f
54p\mod {p^2}\qq\t{if $p=x^2+2y^2\e 1,3\mod 8$,}
\\-\f{45-36(-1)^{\f{p-1}2}}8R_2(p)+\f 32p\mod {p^2}\q\t{if $p\e 5,7\mod {8}$}.
\endcases$$
\endpro

\pro{Conjecture 3.8} Let $p$ be a prime with $p\not=2,3,7$. Then
$$\align &\Ls p3\sum_{k=0}^{p-1}\f{\b{2k}kW_k}{(-27)^k(2k-1)}
\\&\e \cases -\f{2}{3}p+\f{76}{9}y^2\mod {p^2}&\t{if $p=x^2+7y^2\e
1,2,4\mod 7$,}
\\-\f 83\sum_{k=0}^{p-1}\f{\b{2k}k^3}{k+1}-\f{28}9p\mod {p^2}
&\t{if $p\e 3,5,6\mod 7$}
\endcases\endalign$$
and
$$\align &\Ls p3\sum_{k=0}^{p-1}\f{\b{2k}kW_k}{243^k(2k-1)}
\\&\e \cases \f{1676}{81}y^2-\f{142}{81}p\mod {p^2}&\t{if $p=x^2+7y^2\e
1,2,4\mod 7$,}
\\-\f {32}{27}\sum_{k=0}^{p-1}\f{\b{2k}k^3}{k+1}-\f{44}{27}p\mod {p^2}
&\t{if $p\e 3,5,6\mod 7$.}
\endcases\endalign$$
\endpro

\pro{Conjecture 3.9} Let $p$ be a prime with $p\not=2,11$. Then
$$\align &\sum_{k=0}^{p-1}\f{\b{2k}kW_k}{(-44)^k(2k-1)}
\\&\e \cases \f{52}{11}y^2-\f{116}{121}p\mod {p^2}\qq\t{if $p\e 1,3,4,5,9
\mod{11}$ and so $4p=x^2+11y^2$,}
\\\f 9{11}\sum_{k=0}^{p-1}\f{\b{2k}k^2\b{3k}k}{64^k(k+1)}+\f{39}{121}
p\mod {p^2} \q\t{if $p\e 2,6,7,8,10\mod {11}$}
\endcases\endalign$$
\endpro
\pro{Conjecture 3.10} Let $p$ be a prime with $p\not=2,3,19$. Then
$$\align &\Ls p3\sum_{k=0}^{p-1}\f{\b{2k}kW_k}{(-108)^k(2k-1)}
\\&\e \cases -152y^2+\f{4849}{288}p\mod {p^2}&\t{if $\sls p{19}=1$
and so $4p=x^2+19y^2$,}
\\8\ls{-6}p\sum_{k=0}^{p-1}\f{\b{2k}k\b{3k}k\b{6k}{3k}}{(-96)^{3k}(2k-1)}
+\f{241}{288}p\mod {p^2} &\t{if $\ls p{19}=-1$.}
\endcases\endalign$$
\endpro
\pro{Conjecture 3.11} Let $p>3$ be a prime. Then
$$\align&\sum_{k=0}^{p-1}\f{\b{2k}kf_k}{(-4)^k(2k-1)}
\e\cases 2p-4x^2\mod {p^2}&\t{if $p=x^2+3y^2\e 1\mod 3$,}
\\-8R_3(p)\mod {p^2}&\t{if $p\e 2\mod 3$,}
\endcases
\\&\sum_{k=0}^{p-1}\f{\b{2k}kf_k}{50^k(2k-1)} \\&\q\e\cases -\f{13}{25}(4x^2-2p)
-\f{12}{125}(-1)^{\f{p-1}2}p\mod {p^2}&\t{if $p=x^2+3y^2\e 1\mod
3$,}
\\-\f{32}{25}R_3(p)-\f{12}{125}(-1)^{\f{p-1}2}p
\mod {p^2}&\t{if $p\e 2\mod 3$ and $p\not=5$.}
\endcases
\endalign$$\endpro

\pro{Conjecture 3.12} Let $p>3$ be a prime. Then
$$ \align
&\sum_{k=0}^{p-1}\f{\b{2k}kf_k}{96^k(2k-1)}\\& \e\cases
-\f{29}{48}\ls p3(4x^2-2p)-\f p6\mod {p^2}&\t{if $p=x^2+2y^2\e
1,3\mod 8$,}
\\\f{15-12(-1)^{(p-1)/2}}{16}\ls p3R_2(p)+\f p6
\mod {p^2}&\t{if $p\e 5,7\mod 8$.}\endcases\endalign$$
\endpro

\pro{Conjecture 3.13} Let $p$ be a prime with $p\not=2,5$. If
$\sls{-5}p=1$, then
$$\align&\sum_{k=0}^{p-1}\f{\b{2k}kf_k}{16^k(2k-1)}
\\&\e\cases -\f 75x^2+\f{9}{10}p\mod
{p^2}&\t{if $p\e 1,9\mod{20}$ and so $p=x^2+5y^2$,}
\\\f {14}5x^2-\f p2\mod
{p^2}&\t{if $p\e 3,7\mod{20}$ and so $2p=x^2+5y^2$;}\endcases
\endalign$$ if
$\sls{-5}p=-1$, then
$$\sum_{k=0}^{p-1}\f{\b{2k}kf_k}{16^k(2k-1)}
\e
-6(-1)^{\f{p-1}2}\sum_{k=0}^{p-1}\f{\b{2k}k^2\b{4k}{2k}}{(-1024)^k(2k-1)}
-\f{11}8p\mod {p^2}.$$

\pro{Conjecture 3.14} Let $p>3$ be a prime. If $\sls{-6}p=1$, then
 $$\sum_{k=0}^{p-1}\f{\b{2k}kf_k}{32^k(2k-1)}\e
 \cases -\f 74x^2+\f 78p\mod {p^2}
 &\t{if $p=x^2+6y^2\e 1,7\mod
{24}$,}\\-\f 72x^2+\f{7}{8}p\mod {p^2}&\t{if $p=2x^2+3y^2\e 5,11\mod
{24}$;}\endcases$$
 if $\sls{-6}p=-1$, then
$$\sum_{k=0}^{p-1}\f{\b{2k}kf_k}{32^k(2k-1)}\e
\f 94\sum_{k=0}^{p-1}\f{\b{2k}k^2\b{3k}k}{216^k(2k-1)}+\f 14\Ls
p3p\mod {p^2}.$$
\endpro

\pro{Conjecture 3.15} Let $p>3$ be a prime. Then
$$\align&(-1)^{\f{p-1}2}\sum_{k=0}^{p-1}\f{\b{2k}ka_k}{54^k(2k-1)}
\\&\e\cases \f{52}{9}y^2-\f{26+2(-1)^{\f{p-1}2}}{27}p\mod {p^2}&\t{if
$p=x^2+3y^2\e 1\mod 3$,}
\\\f{32}{27}R_3(p)+\f 2{27}(-1)^{\f{p-1}2}p\mod {p^2}&\t{if $p\e 2\mod 3$.}
\endcases\endalign$$
\endpro
\pro{Conjecture 3.16} Let $p>3$ be a prime. Then
$$\sum_{k=0}^{p-1}\f{\b{2k}ka_k}{100^k(2k-1)} \e\cases
-\f{58}{25}x^2+\f{145-18\sls p3}{125}p\mod {p^2}&\t{if $p=x^2+2y^2\e
1,3\mod 8$,}
\\-\f{9}{50}R_2(p)-\f{18}{125}\ls p3p\mod {p^2}&\t{if $p\e 5\mod 8$ and $p\not=5$,}
\\-\f {81}{50}R_2(p)-\f{18}{125}\ls p3p\mod {p^2}&\t{if $p\e 7\mod 8$.}
\endcases$$
\endpro
\pro{Conjecture 3.17} Let $p>3$ be a prime. Then
$$\sum_{k=0}^{p-1}\f{\b{2k}ka_k}{(-12)^k(2k-1)}
\e\cases p-4x^2\mod {p^2}&\t{if $12\mid p-1$ and so $p=x^2+9y^2$,}
\\ 2x^2-p\mod {p^2}&\t{if $12\mid p-5$ and so $2p=x^2+9y^2$,}
\\3\binom{[p/3]}{[p/12]}^2\mod p&\t{if $p\e 7\mod {12}$,}
\\-6\binom{[p/3]}{[p/12]}^2\mod p&\t{if $p\e 11\mod {12}$.}
\endcases$$
\endpro

\pro{Conjecture 3.18} Let $p>3$ be a prime. If $\sls{-6}p=1$, then
 $$\Ls p3\sum_{k=0}^{p-1}\f{\b{2k}ka_k}{36^k(2k-1)}\e
 \cases -\f {14}9x^2+\f 79p\mod {p^2}
 &\t{if $p=x^2+6y^2\e 1,7\mod
{24}$,}\\-\f {28}9x^2+\f{7}{9}p\mod {p^2}&\t{if $p=2x^2+3y^2\e
5,11\mod {24}$;}\endcases$$
 if $\sls{-6}p=-1$, then
$$\sum_{k=0}^{p-1}\f{\b{2k}ka_k}{36^k(2k-1)}\e -2\Ls
p3\sum_{k=0}^{p-1}\f{\b{2k}k^2\b{3k}k}{216^k(2k-1)}-\f 29p\mod
{p^2}.$$
\endpro

\pro{Conjecture 3.19} Let $p$ be an odd prime. Then
$$\sum_{k=0}^{p-1}\f{\b{2k}kG_k}{64^k(k+1)} \e (-1)^{\f{p-1}2}\mod
{p^2} \ \t{and}\ \sum_{k=0}^{p-1}\f{\b{2k}kG_k}{64^k(2k-1)} \e
2(-1)^{\f{p-1}2}p^2\mod {p^3}.
$$\endpro
\pro{Conjecture 3.20} Let $p>3$ be a prime. Then
$$\align &\sum_{k=0}^{p-1}\f{\b{2k}kG_k}{72^k(2k-1)} \e\cases -\f 49x^2+\f
8{27}p\mod {p^2}&\t{if $p=x^2+4y^2\e 1\mod 4$,}
\\-\f 29R_1(p)-\f 2{27}p
\mod {p^2}&\t{if $p\e 3\mod 4$,}
\endcases
\\&\sum_{k=0}^{p-1}\f{\b{2k}kG_k}{576^k(2k-1)}
\e\cases -\f 83x^2+\f {32}{27}p\mod {p^2}&\t{if $p=x^2+4y^2\e 1\mod
4$,}
\\-\f 29R_1(p)+\f 4{27}p
\mod {p^2}&\t{if $p\e 3\mod 4$.}
\endcases
\endalign$$
\endpro

\pro{Conjecture 3.21} Let $p$ be an odd prime. Then
$$\sum_{k=0}^{p-1}\f{\b{2k}kG_k}{128^k(2k-1)} \e \cases -\f
38(4x^2-2p)\mod {p^2}\qq\t{if $p=x^2+2y^2\e 1,3\mod 8$,}
\\-\f {5-4(-1)^{(p-1)/2}}8R_2(p)\mod {p^2}\q\t{if $p\e 5,7\mod 8$.}
\endcases$$
\endpro

\pro{Conjecture 3.22} Let $p>3$ be a prime. Then

$$\align &\sum_{k=0}^{p-1}\f{\b{2k}kG_k}{48^k(2k-1)} \e\cases \f 49x^2\mod
{p^2}&\t{if $p=x^2+3y^2\e 1\mod 3$,}
\\-\f 89R_3(p)-\f 29p\mod {p^2}&\t{if $p\e 2\mod 3$,}
\endcases
\\&\sum_{k=0}^{p-1}\f{\b{2k}kG_k}{(-192)^k(2k-1)}
\e\cases -\f {32}9x^2+\f 43p\mod {p^2}&\t{if $p=x^2+3y^2\e 1\mod
3$,}
\\-\f 89R_3(p)+\f 49p\mod {p^2}&\t{if $p\e 2\mod 3$.}
\endcases
\endalign$$
\endpro

\pro{Conjecture 3.23} Let $p$ be a prime with $p\not=2,3,7$. Then
$$\sum_{k=0}^{p-1}\f{\b{2k}kG_k}{63^k(2k-1)}\e \cases \f
47y^2+\f{26}{1323}p\mod {p^2}\q\t{if $p=x^2+7y^2\e 1,2,4\mod
7$,}\\\f {32}{63}\sum_{k=0}^{p-1}\f{\b{2k}k^3}{k+1}+
\f{1024}{1323}p\mod {p^2}\q\t{if $p\e 3,5,6\mod 7$}
\endcases$$
and
$$\align &\sum_{k=0}^{p-1}\f{\b{2k}kG_k}{(-4032)^k(2k-1)}\\&\q\e \cases \f
{1408}{63}y^2-\f{2368}{1323}p\mod {p^2}\q\t{if $p=x^2+7y^2\e
1,2,4\mod 7$,}\\\f {32}{63}\sum_{k=0}^{p-1}\f{\b{2k}k^3}{k+1}
+\f{688}{1323}p\mod {p^2}\q\t{if $p\e 3,5,6\mod 7$.}\endcases
\endalign$$
\endpro

\pro{Conjecture 3.24} Let $p>3$ be a prime. Then
$$\align&\sum_{k=0}^{p-1}\f{\b{2k}kQ_k}{(-36)^k(2k-1)}
\e\cases -\f 49x^2+\f 29p\mod {p^2}&\t{if $p=x^2+3y^2\e 1\mod 3$,}
\\-\f 89R_3(p)\mod {p^2}&\t{if $p\e 2\mod 3$,}
\endcases
\\&\sum_{k=0}^{p-1}\f{\b{2k}kQ_k}{18^k(2k-1)}
\\&\q\e\cases -\f {52}9x^2+\f {26-12(-1)^{(p-1)/2}}9p\mod {p^2}&\t{if
$p=x^2+3y^2\e 1\mod 3$,}
\\-\f {32}9R_3(p)-\f 43(-1)^{\f{p-1}2}p\mod {p^2}&\t{if $p\e 2\mod 3$.}
\endcases\endalign
$$
\endpro

\pro{Conjecture 3.25} Let $p$ be a prime with $p\not=2,3,11$.
$$\align &\sum_{k=0}^{p-1}\f{\b{2k}kA_k'}{4^k(2k-1)}\\&\e
\cases  2p-2y^2\mod {p^2}\qq\t{if $p\e 1,3,4,5,9\mod {11}$ and so
$4p=x^2+11y^2$,}\\  5\sum_{k=0}^{p-1}\f{\b{2k}k^2\b{3k}{k}}
{64^k(k+1)}+3p\mod {p^2}\q\t{if $p\e 2,6,7,8,10\mod
{11}$.}\endcases\endalign$$
\endpro

\pro{Conjecture 3.26} Let $p$ be a prime with $p\not=2,3,19$. If
$\ls p{19}=1$ and so $4p=x^2+19y^2$, then
$$\sum_{k=0}^{p-1}\f{\b{2k}kA_k'} {36^k(2k-1)} \e
\f{74}{9}y^2-\f{22}{27}p\mod{p^2}.$$ If $\ls p{19}=-1$, then
$$\sum_{k=0}^{p-1}\f{\b{2k}kA_k'} {36^k(2k-1)} \e-\f{40}3
\Ls{-6}p\sum_{k=0}^{p-1}\f{\b{2k}k\b{3k}k\b{6k}{3k}}{(-96)^{3k}(2k-1)}
\mod p.$$
\endpro
\par{\bf Remark 3.1} In [Su3], Z.W. Sun conjectured that for any
prime $p>3$,
$$\align &\sum_{k=0}^{p-1}\f{\b{2k}kA_k'}{4^k}
\e\cases x^2-2p\mod {p^2}&\t{if $\sls p{11}=1$ and so
$4p=x^2+11y^2$,}
\\0\mod {p^2}&\t{if $\sls p{11}=-1$,}\endcases
\\&\sum_{k=0}^{p-1}\f{\b{2k}kA_k'}{36^k}
\e\cases x^2-2p\mod {p^2}&\t{if $\sls p{19}=1$ and so
$4p=x^2+19y^2$,}
\\0\mod {p^2}&\t{if $\sls p{19}=-1$.}\endcases
\endalign$$

\end{document}